\documentclass[onecolumn,amsthm]{autart}   

\usepackage{graphicx}         
\usepackage{amsmath,amssymb,array,cite,xspace}
\usepackage{enumerate,pdfsync,psfrag,bigints}
\usepackage{color}
\usepackage{ifthen}
\usepackage{chngcntr}
\usepackage{psfrag}
\usepackage{marginnote}
\newcommand{\pprime}{{\prime\prime}}

\newcommand{\sigth}[1]{
	\ifthenelse{\equal{#1}{}}{\sigma(\theta)}{\sigma_{#1}(\theta)}
}
\newcommand{\dsigth}[1]{
	\ifthenelse{\equal{#1}{}}{\sigma^{\prime}(\theta)}{\sigma^{\prime}_{#1}(\theta)}
}
\newcommand{\ddsigth}[1]{
	\ifthenelse{\equal{#1}{}}{\sigma^{\pprime}(\theta)}{\sigma^{\pprime}_{#1}(\theta)}
}

\newcounter{margincomments}\stepcounter{margincomments}

\newcommand{\cK}{\mathcal K}

\newcommand{\cH}{\mathcal H}
\newcommand{\cX}{\mathcal X}
\newcommand{\cQ}{\mathcal Q}
\newcommand{\cU}{\mathcal U}
\newcommand{\cV}{\mathcal V}
\newcommand{\cW}{\mathcal W}

\newcommand{\T}{{\mathscr T}}
\renewcommand{\Re}{\mathbb R}
\newcommand{\Se}{\mathbb S}

\newcommand{\bperp}{B^\perp}
\newcommand{\trans}{^\top}

\newcommand{\vhc}{VHC\xspace}
\newcommand{\vhcs}{VHCs\xspace}

\newcommand{\ucirc}{\mathbb{S}^1}

\newcommand{\cI}{\mathcal{I}}

\DeclareMathOperator{\rank}{rank}
\DeclareMathOperator{\image}{Im}

\DeclareMathOperator{\diag}{diag}

\DeclareMathOperator{\col}{col}

\DeclareMathOperator{\Hess}{Hess}

\newcommand{\z}{z}
\newcommand{\s}{\vartheta}
\renewcommand{\T}{{\tilde T}}
\newcommand{\Tconstraint}{{T_1}}
\newcommand{\Tgamma}{{T_2}}

\newenvironment{definition}{
	\begin{defn}}{\hfill $\bigtriangleup$
	\end{defn}
}
\begin{document}
\begin{frontmatter}

\title{Dynamic Virtual Holonomic Constraints for Stabilization of Closed Orbits in Underactuated Mechanical Systems}

\thanks[footnoteinfo]{This paper was not presented at any IFAC 
meeting. Corresponding author M.~Maggiore. Tel. (Fax) +1-416-946-5095.}

\author[UToronto]{Alireza Mohammadi}\ead{alireza.mohammadi@mail.utoronto.ca},    
\author[UToronto]{Manfredi Maggiore}\ead{maggiore@control.utoronto.ca},               
\author[UParma]{Luca Consolini}\ead{lucac@ce.unipr.it}  

\address[UToronto]{Department of Electrical and Computer Engineering, University of Toronto, 10 King's College Road, Toronto,
	Ontario, M5S 3G4, Canada.}  
\address[UParma]{Dipartimento di Ingegneria dell'Informazione, Via Usberti 181/a, 43124 Parma, Italy}             

\begin{keyword}                          
Dynamic virtual holonomic constraints; underactuated mechanical systems; orbital stabilization.               
\end{keyword}

\begin{abstract}             
This article investigates the problem of enforcing a virtual holonomic
constraint (\vhc) on a mechanical system with degree of underactuation
one while simultaneously stabilizing a closed orbit on the constraint
manifold. This problem, which to date is open, arises when designing
controllers to induce complex repetitive motions in robots.  In this
paper, we propose a solution which relies on the parameterization of
the \vhc by the output of a double integrator. While the original
control inputs are used to enforce the \vhc, the control input of the
double-integrator is designed to asymptotically stabilize the closed
orbit and make the state of the double-integrator converge to
zero. The proposed design is applied to the problem of making a PVTOL
aircraft follow a circle on the vertical plane with a desired speed
profile, while guaranteeing that the aircraft does not roll over
for suitable initial conditions.
\end{abstract}
\end{frontmatter}

\nocite{FreLaHMetRobShiJoh09}

Virtual holonomic constraints (\vhcs) have been recognized to be key
to solving complex motion control problems in robotics. There is an
increasing body of evidence from bipedal
robotics~\cite{grizzle2001asymptotically,grizzle2014models,
  WesGriCheChoMor07}, snake robot
locomotion~\cite{mohammadi2015maneuver}, and repetitive motion
planning~\cite{ahmed2013kite,shiriaev2007virtual} that \vhcs
constitute a new motion control paradigm, an alternative to the
traditional reference tracking framework. The key difference with the
standard motion control paradigm of robotics is that, in the \vhc
framework, the desired motion is parameterized by the states of the
mechanical system, rather than by time.

Geometrically, a \vhc is a subset\footnote{More precisely, an
  embedded submanifold.} of the configuration manifold of the
mechanical system. Enforcing a \vhc means stabilizing the subset of
the state space of the mechanical system where the generalized
coordinates of the mechanical system satisfy the \vhc, while the
generalized velocity is tangent to the \vhc. This subset is called
the {\em constraint manifold}.  

Grizzle and collaborators (see, e.g.,~\cite{WesGriCheChoMor07}) have
shown that the enforcement of certain \vhcs on a biped robot leads,
under certain conditions, to the orbital stabilization of a hybrid
closed orbit corresponding to a repetitive walking gait. The orbit in
question lies on the constraint manifold, and the mechanism
stabilizing it is the dissipation of energy that occurs when a foot
impacts the ground. In a mechanical system without impacts, this
stabilization mechanism disappears, and the enforcement of the \vhc
alone is insufficient to achieve the ultimate objective of stabilizing
a repetitive motion. Some researchers~\cite{ShiPerWit05,ShiFreGus10}
have addressed this problem by using the \vhc exclusively for motion
planning, i.e., to find a desired closed orbit. Once a suitable closed
orbit is found, a time-varying controller is designed by linearizing
the control system along the orbit. In this approach, the constraint
manifold is not an invariant set for the closed-loop system, and thus
the \vhc is not enforced via feedback.

To the best of our knowledge, for mechanical control systems with
degree of underactuation one, the problem of simultaneous enforcement
of a \vhc and orbital stabilization of a closed orbit lying on the
constraint manifold is still open. The challenge in addressing this
problem lies in the fact that the dynamics of the mechanical control
system on the constraint manifold are unforced. Therefore, any
feedback that asymptotically stabilizes the desired closed orbit
cannot render the constraint manifold invariant, and thus cannot
enforce the \vhc. To overcome this difficulty, in this paper we
propose to render the \vhc{} {\em dynamic}. By doing that, under
suitable assumptions it is possible to stabilize the desired closed
orbit while simultaneously enforcing the dynamic \vhc.

\textbf{Contributions of the paper.} This paper presents the first
solution of the simultaneous stabilization problem just described for
mechanical control systems with degree of underactuation one.
Leveraging recent results in~\cite{MohMagCon17}, we consider \vhcs
that induce Lagrangian constrained dynamics. The closed orbits on the
constraint manifold are level sets of a ``virtual'' energy
function. We make the \vhc dynamic by parametrizing it by the output
of a double-integrator. We use the original control inputs of the
mechanical system to stabilize the constraint manifold associated with
the dynamic \vhc, and we use the double-integrator input to
asymptotically stabilize the selected orbit on the constraint
manifold. Because the output of the double-integrator acts as a
perturbation of the original constraint manifold, we also make sure
that the state of the double-integrator converges to zero. To achieve
these objectives, we develop a novel theoretical result giving
necessary and sufficient conditions for the exponential
stabilizability of closed orbits for control-affine systems.

The benefits associated with the simultaneous stabilization
proposed in this paper are as follows. First, in the proposed
framework one may assign the speed of convergence of solutions to
the constraint manifold independently of the orbit stabilization
mechanism. In particular, one may enforce the dynamic \vhc
arbitrarily fast\footnote{Naturally, actuator saturation will limit
  the maximum attainable speed of convergence to the constraint
  manifold.}, so that after a short transient, the qualitative
behaviour of trajectories of the closed-loop system is determined by
the dynamic \vhc.  Second, since the constraint manifold is
asymptotically stable for the closed-loop system, trajectories
originating near the constraint manifold remain close to it
thereafter. From a practical standpoint, the two features just
highlighted mean that the dynamic \vhc offers some control over the
transient behaviour of the closed-loop system.  The simultaneous
stabilization of the closed orbit means that, without violating the
dynamic \vhc, an extra stabilization mechanism makes the trajectories
of the closed-loop system converge to the closed orbit.

The property just described is illustrated in this paper with an
example, the model of a PVTOL aircraft moving along a unit circle on
the vertical plane. The control specification is to make the aircraft
traverse the circle with bounded speed, while guaranteeing that the
aircraft does not undergo full revolutions along its longitudinal
axis.  In this context, the \vhc constrains the roll angle of the
aircraft as a function of its position on the circle, preventing the
aircraft from rolling over. On the other hand, the simultaneous
stabilization of the closed orbit corresponds to stabilizing a desired
periodic speed profile on the circle {\em without violating the
  constraint.} The double-integrator state perturbs the constraint so
as to induce the orbit stabilization mechanism.

\textbf{Relevant literature.}  Previous work employs \vhcs to
stabilize desired closed orbits for underactuated mechanical
systems~\cite{canudas2002orbital,canudas2004concept,ShiPerWit05,FreLaHMetRobShiJoh09}. Canudas-de-Wit
and collaborators~\cite{canudas2002orbital} propose a technique to
stabilize a desired closed orbit that relies on enforcing a virtual
constraint and on dynamically changing its geometry so as to impose
that the reduced dynamics on the constraint manifold match the
dynamics of a nonlinear
oscillator. In~\cite{canudas2004concept,ShiPerWit05}, Canudas-de-Wit,
Shiriaev, and collaborators employ \vhcs to aid the selection of
closed orbits of underactuated mechanical systems. It is demonstrated
that an unforced second-order system possessing an integral of motion
describes the constrained motion. Assuming that this unforced system
has a closed orbit, a linear time-varying controller is designed that
yields exponential stability of the closed orbit. With the exception
of~\cite{canudas2002orbital}, the papers above do not guarantee the
invariance of the \vhc for the closed loop system.  The idea of
event-triggered dynamic \vhcs has appeared in the work by Morris and
Grizzle in~\cite{morris2009hybrid} where the authors construct a
hybrid invariant manifold for the closed-loop dynamics of biped robots
by updating the \vhc parameters after each impact with the ground.
This approach is similar in spirit to the one presented in this paper.
Finally, the paper~\cite{CelAnd16} discusses collocated VHCs, i.e.,
VHCs parametrized by actuated variables.  In
Section~\ref{sec:stab:shir}, we discuss the differences between the
method presented in this article and the ones
in~\cite{canudas2002orbital,canudas2004concept,ShiPerWit05,ShiFreGus10}. We
also discuss the conceptual similarities between the method presented
in this article and the one in~\cite{morris2009hybrid}.


\textbf{Organization.} This article is organized as follows.  We
review preliminaries in Section~\ref{sec:stab:prelim}. The formal
problem statement and our solution strategy are presented in
Section~\ref{sec:stab:prob_stat}.  In Section~\ref{sec:stab:dynConstr}
we present dynamic \vhcs. In Section~\ref{sec:stab:s_control} we
present a novel result of a general nature providing necessary and
sufficient conditions for the exponential stabilizability of closed
orbits for control-affine systems, and use it to design the input of
the double-integrator to stabilize the closed orbit relative to the
constraint manifold. In Section~\ref{sec:stab:sol} we present the
complete control law solving the \vhc-based orbital stabilization
problem.  In Section~\ref{sec:stab:shir} we discuss the differences
between the method presented in this article and the ones
in~\cite{canudas2004concept,canudas2002orbital,ShiPerWit05}. Finally,
in Section~\ref{sec:MotEx} we apply the ideas of this paper to a path
following problem for the PVTOL aircraft.

\textbf{Notation.} If $x \in \Re$ and $T>0$, then $x$ modulo $T$ is
denoted by $[x]_T$, and the set $\{[x]_T: x\in \Re\}$ is denoted by
$[\Re]_T$. This set can be given a manifold structure which makes it
diffeomorphic to the unit circle $\ucirc$. If $a$ and $b$ are vectors,
then $\col(a,b):=[a\trans \ b\trans]\trans$. If $a,b\in \Re^n$, we
denote $\langle a,b \rangle = a\trans b$, and $\|a\| = \langle a,a
\rangle^{1/2}$. If $A \in \Re^{n\times n}$, we denote by $\|A\|_2$ the
induced two-norm of $A$. If $(\cX,d)$ is a metric space, $\Gamma$ is a
subset of $\cX$, and $x \in \cX$, we denote by $\|x\|_\Gamma$ the
point-to-set distance of $x$ to $\Gamma$, defined as $\|x\|_\Gamma
:=\inf_{y \in \Gamma} d(x,y)$.

If $h: M \to N$ is a smooth map between smooth manifolds, and $q \in
M$, we denote by $dh_q: T_q M \to T_{h(q)} N$ the derivative of $h$ at
$q$ (in coordinates, this is the Jacobian matrix of $h$ evaluated at
$q$), and if $M$ has dimension 1, then we may use the notation $h'(q)$
in place of $dh_q$.  If $M_1,M_2,N$ are smooth manifolds and $f:M_1
\times M_2 \to N$ is a smooth function, then $\partial_{q_1}
f(q_1,q_2)$ denotes the derivative of the map $q_1 \mapsto f(q_1,q_2)$
at $q_1$. If $f: M \to TM$ is a vector field on $M$ and $h:M \to
\Re^m$ is $C^1$, then $L_f h: M \to \Re^m$ is defined as $L_f
h(q):=dh_q f(q)$.  For a function $h: M \to \Re^m$, we denote by
$h^{-1}(0) :=\{q \in M: h(q)=0\}$.

If $A\in \Re^{m \times n}$ has full row-rank, we denote by
$A^{\dagger}$ the pseudoinverse of $A$, $A^{\dagger}=A^{\top}(A
A^{\top})^{-1}$.  Given a $C^2$ scalar function $f:\Re^n \to \Re$, we
denote by $\mbox{Hess}(f)$ its Hessian matrix.


%

%
\section{Preliminaries}\label{sec:stab:prelim}

Consider the underactuated mechanical control system 
\begin{equation}\label{eq:stab:ELsys}
D(q) \ddot q +  C(q,\dot q) \dot q + \nabla P(q) = B(q) \tau,
\end{equation}
where $q=\big(q_1,\ldots,q_n\big)\in\cQ$ is the configuration vector
with $q_i$ either a displacement in $\Re$ or an angular variable in
$[\Re]_{T_i}$, with $T_i>0$. The configuration space $\cQ$ is,
therefore, a generalized cylinder. In \eqref{eq:stab:ELsys}, $B : \cQ
\to \Re^{n \times n-1}$ is $C^1$ and it has full rank $n-1$.  Also,
$D(q)$, the inertia matrix, is positive definite for all $q$, and
$P(q)$, the potential energy function, is $C^1$. We assume that there
exists a left-annihilator of $B(q)$; specifically, there is a $C^1$
function $\bperp:\cQ \to \Re^{1\times n} \backslash \{0\}$ such that
$\bperp(q) B(q) =0$ for all $q\in \cQ$.
\begin{definition}[\hspace*{-4pt}\cite{Maggiore-2013VHC}]\label{defn:vhc} 
A relation $h(q)=0$, where $h: \cQ \to \Re^{k}$ is $C^2$, is a {\bf
  regular virtual holonomic constraint (\vhc) of order $k$}
for system~\eqref{eq:stab:ELsys}, if~\eqref{eq:stab:ELsys} with output
$e=h(q)$ has well-defined vector relative degree $\{2,\cdots,2\}$
everywhere on the {\bf constraint manifold}
\begin{equation}\label{eq:constraint_manifold}
\Gamma := \{ (q,\dot q): h(q) =0, \ dh_q \dot q=0\},
\end{equation}
i.e., the matrix $dh_q D^{-1}(q) B(q)$ has full row rank for all $q
\in h^{-1}(0)$.
\end{definition}
The constraint manifold $\Gamma$ in~\eqref{eq:constraint_manifold} is
just the zero dynamics manifold associated with the output
$e=h(q)$. For a \vhc of order $n-1$, the set $h^{-1}(0)$ is a
collection of disconnected regular curves, each one diffeomorphic to
either the unit circle or the real line. From now on, we will assume
that $h^{-1}(0)$ is diffeomorphic to $\Se^1$. 

Necessary and sufficient conditions for a
relation $h(q)=0$ to be a regular \vhc of order $n-1$ are given in the
following proposition.
\begin{prop}[\hspace*{-4pt}\cite{Maggiore-2013VHC}]
\label{prop:regular_vhc}
Let $h: \cQ \to \Re^{n-1}$ be $C^2$ and such that $\rank dh_q = n-1$
for all $q \in h^{-1}(0)$. Then, $h(q)=0$ is a regular \vhc of order $n-1$ for
system~\eqref{eq:stab:ELsys} if and only if for each $q \in
h^{-1}(0)$,
\[
T_q h^{-1}(0) \oplus  \image (D^{-1}(q)B(q)) = T_q \cQ.
\]
Moreover, if $\sigma: \Re
\rightarrow \cQ$ is a regular parameterization of $h^{-1}(0)$, then
$h(q)=0$ is a regular \vhc for system~\eqref{eq:stab:ELsys} if and only
if
\[
(\forall \theta \in \Re) \ \bperp(\sigth{}) D(\sigth{}) \dsigth{} \neq
0.
\]
\end{prop}
By definition, if $h:\cQ \to \Re^{n-1}$ is a regular \vhc,
system~\eqref{eq:stab:ELsys} with output $e=h(q)$ has vector relative
degree $\{2,\ldots,2\}$. In order to asymptotically stabilize the
constraint manifold, one may employ an input-output linearizing
feedback, as detailed in the next proposition.  Before stating the
proposition, we define the notion of set stability used in this paper.
\begin{definition}
  Consider a dynamical system $\Sigma$ on a metric space $(\cX,d)$ with
  continuous local flow map $\phi(t,x_0)$, defined on an open subset of
  $\Re \times \cX$. The set $\Gamma \subset \cX$ is {\bf stable} for
  $\Sigma$ if for all $\varepsilon >0$ there exists a neighbourhood $U$
  of $\Gamma$ such that for all $x_0 \in U$ such that $\phi(t,x_0)$ is
  defined for all $t \geq 0$, $\| \phi(t,x_0) \|_\Gamma < \varepsilon$
  for all $t \geq 0$. The set $\Gamma$ is {\bf asymptotically stable}
  for $\Sigma$ if it is stable and there exists a neighbourhood $U$ of
  $\Gamma$ such that for all $x_0 \in U$ such that $\phi(t,x_0)$ is
  defined for all $t \geq 0$, $\|\phi(t,x_0)\|_\Gamma \to 0$ as $t \to
  \infty$.
\end{definition}
We remark that if $\gamma \subset \cX$ is a closed orbit of the
dynamical system $\Sigma$, then the notion of asymptotic stability of
$\gamma$ coincides with that of {\em asymptotic orbital stability}
found in the literature (see,
e.g.,~\cite[Definition~8.2]{Khalil-2002}). Therefore, in the sequel we
will speak of asymptotic stability of a closed orbit $\gamma$.
\begin{prop}[\hspace*{-4pt}\cite{Maggiore-2013VHC}]\label{prop:stabilization}
Let $h(q)=0$ be a regular \vhc of order $n-1$ for
system~\eqref{eq:stab:ELsys} with associated constraint manifold
$\Gamma$ in~\eqref{eq:constraint_manifold}.  Let $H(q,\dot q) =
\col(h(q),dh_q \dot q)$, and assume that there exist two class-$\cK$
functions $\alpha_1,\alpha_2$ such that
\begin{equation}\label{eq:classK}
\alpha_1(\|(q,\dot q)\|_\Gamma) \leq H(q,\dot q) \leq
\alpha_2(\|(q,\dot q)\|_\Gamma).
\end{equation}
Let $A(q)=dh_q D^{-1}(q) B(q)$, $e=h(q)$, and $\cH =
\col(\cH_1,\ldots, \cH_{n-1})$, where $\cH_i = \dot q\trans \Hess(h_i)
\dot q$.  Then, for all $k_p,k_d>0$, the input-output linearizing
controller
\[
\tau = A^{-1}(q) \left\{ dh_q D^{-1}(q) [C(q,\dot q) \dot q + \nabla
  P(q)] -\cH(q,\dot q) - k_p e - k_d \dot e\right\},
\]
asymptotically stabilizes the constraint manifold $\Gamma$.
\end{prop}
Once the constraint manifold $\Gamma$ has been rendered invariant by
the above feedback, the motion on $\Gamma$ is described by a
second-order unforced differential equation, as detailed in the next
proposition.
\begin{prop}[\hspace*{-1ex}\cite{WesGriKod03,ShiPerWit05,MohMagCon17}]
\label{prop:red_dyn}
Let $h(q)=0$ be a regular \vhc of order $n-1$ for
system~\eqref{eq:stab:ELsys}. Assume that $h^{-1}(0)$ is diffeomorphic
to $\Se^1$. For some $\Tconstraint>0$, let $\sigma: [\Re]_{\Tconstraint} \rightarrow
\cQ$ be a regular parameterization of $h^{-1}(0)$. Then letting
$(q,\dot q) = ( \sigma(\theta), \sigma'(\theta)\dot \theta)$, the
dynamics on the set $\Gamma$ in~\eqref{eq:constraint_manifold} are
globally described by
\begin{equation}\label{eq:constrained_dynamics}
\ddot \theta = \Psi_1(\theta) + \Psi_2(\theta) \dot{\theta}^2,
\end{equation}
where $(\theta,\dot \theta) \in [\Re]_\Tconstraint \times \Re$ and
\begin{equation}\label{eq:Psi_functions}
\begin{aligned}
& \Psi_1(\theta) = -\frac{\bperp\nabla P}{\bperp D \sigma^{\prime}}\bigg|_{q=\sigth{}},  \\
& \Psi_2(\theta) = -\frac{\bperp D \sigma^{\pprime} + \sum_{i=1}^n
    \bperp_i \sigma^{\prime\trans} Q_i \sigma^{\prime}}{\bperp D
    \sigma^{\prime}}\Bigg|_{q=\sigth{}},
\end{aligned}
\end{equation}
and where $\bperp_i$ is the $i$-th component of $\bperp$ and $Q_i$ is
an $n \times n$ matrix whose $(j,k)$-th component is $(Q_i)_{jk}
=(1/2) (\partial_{q_k} D_{ij} + \partial_{q_j} D_{ik} - \partial_{q_i}
D_{kj})$.
\end{prop}

Henceforth, we will refer to~\eqref{eq:constrained_dynamics} as the
{\bf reduced dynamics}. System~\eqref{eq:constrained_dynamics} is
unforced since all $n-1$ control directions are used to make the
constraint manifold $\Gamma$ invariant. In the context of nonlinear
control, the reduced dynamics~\eqref{eq:constrained_dynamics} are a
coordinate representation of the zero dynamics vector field of the
mechanical system~\eqref{eq:stab:ELsys} with output $e = h(q)$.
Proposition~\ref{prop:red_dyn} is a direct consequence of the fact
that outputs of a nonlinear systems with a well-defined relative
degree induce a globally-defined zero dynamics vector
field~\cite{Isi95}.

\begin{figure}[htb]
\label{fig:constraint_manifold}
\psfrag{t}{$\theta$}
\psfrag{T}{$\dot \theta$}
\psfrag{G}{$\Gamma$}
\psfrag{S}{$q=\sigma(\theta)$}
\psfrag{x}{$q_1$}
\psfrag{y}{$q_n$}
\psfrag{h}[c]{$h^{-1}(0)$}
\centerline{\includegraphics[height=.2\textheight]{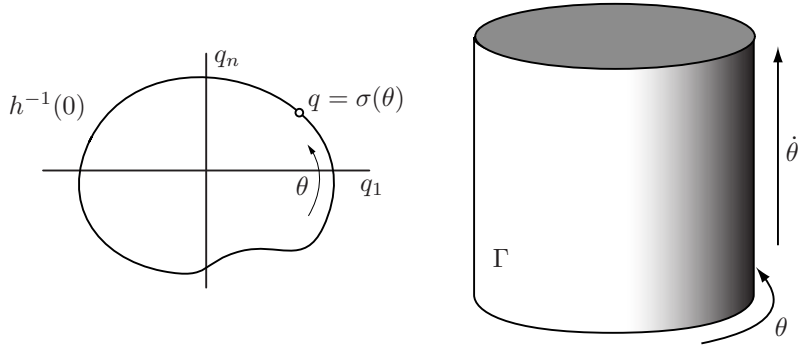}}
\caption{A \vhc $h(q)=0$ and the associated constraint manifold
  $\Gamma$. The \vhc, on the left-hand side, is a curve diffeomorphic
  to $\Se^1$, with parametrization $q=\sigma(\theta)$. The associated
  constraint manifold $\Gamma$, on the right-hand side, is a cylinder
  parametrized by the variables $(\theta,\dot \theta)$ via the
  diffeomorphism $(q,\dot q) = (\sigma(\theta),\sigma'(\theta)\dot
  \theta)$.}
\end{figure}

\begin{rem} The proposition above states that the second-order
  differential equation~\eqref{eq:constrained_dynamics} with state
  space $[\Re]_{\Tconstraint} \times \Re$ represents the dynamics on
  the set $\Gamma$. The geometric underpinning of this statement is
  the fact that $\Gamma$ is diffeomorphic to the cylinder
  $[\Re]_{\Tconstraint} \times \Re$ via the diffeomorphism
  $[\Re]_{\Tconstraint} \times \Re \to \Gamma$, $(\theta,\dot
  \theta) \mapsto (\sigma(\theta), \sigma'(\theta) \dot \theta)$. We
  can therefore identify $\Gamma$ with the cylinder
  $[\Re]_{\Tconstraint} \times \Re$, and parametrize it with the
  variables $(\theta,\dot \theta)$, see
  Figure~\ref{fig:constraint_manifold}.\hfill
$\triangle$
\label{rem:cylinder}\end{rem}

Under certain conditions, the reduced
dynamics~\eqref{eq:constrained_dynamics} have a Lagrangian
structure. Define
\begin{equation}
\label{eq:stab:M_V}
M(\theta) := \exp\left(-2\int_0^{\theta}\Psi_2(\tau)d\tau\right),\;
V(\theta):= -\int_0^{\theta} \Psi_1(\tau) M(\tau)d\tau.
\end{equation}
\begin{prop}[\hspace*{-4pt} \cite{MohMagCon17}]  
\label{prop:ESAIM} 
Consider the reduced dynamics~\eqref{eq:constrained_dynamics} with
state space $[\Re]_{\Tconstraint} \times
\Re$. System~\eqref{eq:constrained_dynamics} is Lagrangian if and only
if the functions $M(\cdot)$ and $V(\cdot)$ in~\eqref{eq:stab:M_V} are
$\Tconstraint$-periodic, in which case the Lagrangian function is
given by
$L(\theta,\dot{\theta})=(1/2)M(\theta)\dot{\theta}^2-V(\theta)$.
\end{prop}

An immediate consequence of the foregoing result is that, when the
reduced dynamics~\eqref{eq:constrained_dynamics} are Lagrangian, the
orbits of~\eqref{eq:constrained_dynamics} are characterized by the
level sets of the energy function
\begin{equation}\label{eq:energy}
E(\theta,\dot{\theta}) = \frac 1 2 M(\theta) \dot{\theta}^2 +
V(\theta).
\end{equation}
We remark that the energy function $E(\theta,\dot \theta)$ appeared in
the work~\cite{FreMetShiSpo09}. A different function, dependent on
initial conditions, was presented in~\cite{ShiPerWit05} as an
``integral of motion'' of the reduced
dynamics~\eqref{eq:constrained_dynamics}.

Almost all orbits of the reduced
dynamics~\eqref{eq:constrained_dynamics} are closed, and they belong
to two distinct families, defined next.

\begin{figure}[htb]
\label{fig:closed_orbits}
\psfrag{t}{$\theta$}
\psfrag{T}{$\dot \theta$}
\psfrag{a}{$\gamma_1$}
\psfrag{G}{$\Gamma$}
\psfrag{b}{$\gamma_2$}
\psfrag{x}{$q_1$}
\psfrag{y}{$q_n$}
\centerline{\includegraphics[height=.15\textheight]{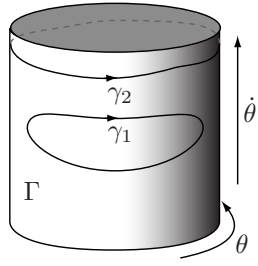}}
\caption{An illustration of the two types of closed orbits on
    $\Gamma$ exhibited by the reduced
  dynamics~\eqref{eq:constrained_dynamics} under the assumptions of
  Proposition~\ref{prop:ESAIM}. The orbit $\gamma_1$ is an
  oscillation, while $\gamma_2$ is a rotation.}
\end{figure}

\begin{definition}
A closed orbit $\gamma$ of the reduced
dynamics~\eqref{eq:constrained_dynamics} is said to be a {\bf rotation
  of} $\theta$ if $\gamma$ is homeomorphic to a circle $\{(\theta,
\dot{\theta}) \in [\Re]_{T} \times \Re:
\dot{\theta}=\text{constant}\}$ via a homeomorphism of the form
$(\theta, \dot{\theta}) \mapsto (\theta,T(\theta)\dot{\theta})$;
$\gamma$ is an {\bf oscillation of} $\theta$ if it is homeomorphic to
a circle $\{(\theta,\dot{\theta}) \in [\Re]_{\Tconstraint} \times \Re:
\theta^2 + \dot{\theta}^2 = \text{constant}\}$ via a homeomorphism of
the form $(\theta, \dot{\theta}) \mapsto
(\theta,T(\theta)\dot{\theta})$.
\end{definition}
In~\cite[Proposition~4.7]{Maggiore-2013VHC}, it is shown that if the
assumptions of Proposition~\ref{prop:red_dyn} hold, then almost all
orbits of~\eqref{eq:constrained_dynamics} are either oscillations or
rotations. Oscillations and rotations are illustrated in
Figure~\ref{fig:closed_orbits}. It is possible to give an explicit
regular parameterization of rotations and oscillations which will be
useful in what follows.

If $\gamma$ is a rotation with associated energy value $E_0$, then we
may solve $E(\theta,\dot \theta)=E_0$ for $\dot \theta$ obtaining
\[
\dot \theta = \pm \sqrt{\frac{2}{M(\theta)} (E_0 - V(\theta))},
\]
with plus sign for counterclockwise rotation, and minus sign for
clockwise rotation. Thus a rotation $\gamma$ is the graph of a
function, which leads to the natural regular parameterization
$[\Re]_{\Tconstraint} \to [\Re]_{\Tconstraint \times \Re}$ given by
\begin{equation}\label{eq:param_rotation}
\s \mapsto (\varphi_1(\s),\varphi_2(\s)) = \left(\s, \pm
\sqrt{\frac{2}{M(\s)} (E_0 - V(\s))}\right).
\end{equation}
Concerning oscillations, it was shown
in~\cite[Lemma~3.12]{consolini2010path} that they are mapped
homeomorphically to circles via the homeomorphism
\[
(\theta,\dot \theta) \mapsto \left( \theta, 
T(\theta)\dot \theta\right), \quad T(\theta) = \sqrt{ \frac{R^2 - (\theta -
    C)^2}{\frac{2}{M(\theta)}(E_0 - V(\theta))}}.
\]
In the above, $E_0$ is the energy level associated with $\gamma$, and
\[
\theta^1 := \min_{(\theta,\dot \theta) \in \gamma}\theta, \ \theta^2
:= \max_{(\theta,\dot \theta) \in \gamma}\theta,
\ C:=(\theta^1+\theta^2)/2, \ R:=(\theta^2 - \theta^1)/2.
\]
The image of $\gamma$ under the above homeomorphism is a circle of
radius $R$ centred at $(C,0)$. Using this fact, we get the
following regular parameterization $[\Re]_{2\pi} \to
[\Re]_{\Tconstraint} \times \Re$:
\begin{equation}
\label{eq:param_oscillation}
\s \mapsto ( \varphi_1(\s),\varphi_2(\s) ) =\left(C+R \cos(\s),
\frac{R \sin(\s)}{T(C+R\cos(\s))} \right).
\end{equation}

\section{Problem formulation}\label{sec:stab:prob_stat}

Consider the mechanical control system~\eqref{eq:stab:ELsys} with $n$
DOFs and $n-1$ control inputs. Let $h(q)=0$ be a regular \vhc of order
$n-1$, and assume that $h^{-1}(0)$ is diffeomorphic to $\Se^1$. As
before, let $\sigma:[\Re]_{\Tconstraint}\to \cQ$, $\Tconstraint>0$, be
a regular parameterization of $h^{-1}(0)$.

Assume that the dynamics~\eqref{eq:constrained_dynamics} are
Lagrangian, so that almost all of its closed orbits are rotations or
oscillations. In particular, almost every orbit of the reduced
dynamics on the constraint manifold is closed, and it corresponds to a
certain speed profile.  Pick one such orbit of interest\footnote{Here
  we assume that the level set $\{(\theta,\dot \theta): E(\theta,\dot
  \theta) = E_0\}$ is connected. There is no loss of generality in
  this assumption, since the theory developed below relies on the
  regular
  parametrizations~\eqref{eq:param_rotation},~\eqref{eq:param_oscillation},
  which one can use to select one of the desired connected components
  of $\gamma$.}, $\gamma=\{(\theta,\dot \theta): E(\theta,\dot \theta)
= E_0\}$. As pointed out in the introduction, since the reduced
dynamics~\eqref{eq:constrained_dynamics} are unforced, it is
impossible to stabilize this orbit while preserving the invariance of
the constraint manifold. 

The idea we explore in this paper is to introduce a dynamic
perturbation of the constraint manifold.  We define a one-parameter
family of \vhcs $h^s(q)=0$, where $s \in \Re$ is the parameter and the
map $h^s(q)$ is such that $h^0(q) = h(q)$.  In
Section~\ref{sec:stab:dynConstr} of this paper we choose the
parametrization $h^s(q):=h(q - Ls)$, where $L \in \Re^n$ is a
parameter vector, but other parametrizations are possible. Each
\vhc in the one-parameter family can be viewed as a perturbation of
the original \vhc $h(q)=0$.  The parameter $s$ is dynamically
adapted by means of the double-integrator $\ddot s=v$, where $v$ is
a new control input.  System~\eqref{eq:stab:ELsys}, augmented with
this double-integrator, has state $(q,\dot q,s,\dot s) \in \cQ \times
\Re^n \times \Re \times \Re$, and control input $(\tau,v) \in
\Re^{n-1} \times \Re$.  We will employ $\tau$ to stabilize the
constraint manifold associated with the dynamic \vhc $h^s(q)=0$, and
$v$ to stabilize a new closed orbit for the augmented system. We
will detail our solution steps below, but first we will formulate
precisely our control specifications.

We begin by defining the constraint manifold associated with the
family of \vhcs $h^s(q)=0$ as
\[
\bar \Gamma = \{(q,\dot q,s,\dot s) : h^s(q)=0, \partial_q h^s \dot q
+ \partial_s h^s \dot s =0\}.
\]
The original constraint manifold $\Gamma$ is embedded in $\bar \Gamma$
as the intersection of $\bar \Gamma$ with the plane $\big\{
(q,\dot{q}, s , \dot{s}) :\, s = 0 ,\, \dot{s} = 0 \big\}$, because
the identity $h^0(q) =h(q)$ implies that
$
\{(q,\dot q,s,\dot s): (q,\dot q) \in \Gamma, (s,\dot s)=(0,0)\}
\subset \bar \Gamma$.
For the one-parameter family $h^s(q) = h(q -Ls)$ used in this paper,
we will show in the proof of
Proposition~\ref{prop:stab:stabilizability} that the distance of a
point $(q,\dot q,s,\dot s) \in \cQ \times \Re^n \times \Re \times \Re$
to the set $\bar \Gamma$ can be expressed as
\[
\| (q,\dot q,s,\dot s)\|_{\bar \Gamma} = \|(q-Ls,\dot q - L \dot
s)\|_\Gamma,
\]
from which one readily deduces the inequality
\[
\|(q,\dot q)\|_{\Gamma} \leq \| (q,\dot q,s,\dot s)\|_{\bar \Gamma}+
\|(Ls,L\dot s)\|.
\]
Thus, if the state $(q,\dot q,s,\dot s)$ of the augmented system is
close to $\bar \Gamma$, and if $\|(s,\dot s)\|$ is small, the state
$(q,\dot q)$ of the mechanical system is close to $\Gamma$, the
original constraint manifold. For this reason, our first control
specification is the asymptotic stabilization of $\bar \Gamma$.

The second control specification for the augmented system will
correspond, in an appropriate manner, to the  orbital stabilization  of
$\gamma$. The curve $\gamma$ is contained in the state space of the
original mechanical system, so we need to lift it to the state space
of the augmented system. The lift in question is
\[
\bar \gamma := \{(q,\dot q,s,\dot s) : s=\dot s=0, \, q =
\sigma(\theta), \, \dot q = \sigma'(\theta) \dot \theta, \,
(\theta,\dot \theta) \in [\Re]_{\Tconstraint} \times \Re, \,
E(\theta,\dot \theta)=E_0\}.
\]
The second control specification is the asymptotic stabilization of
$\bar \gamma$.  That this is indeed the right control specification
follows from the observation that if $(q,\dot q,s,\dot s) \in \bar
\gamma$, then $(q,\dot q)\in \gamma$. Thus, the stabilization of $\bar
\gamma$ for the augmented system makes the trajectories of the
mechanical system~\eqref{eq:stab:ELsys} converge to $\gamma$, as
desired. Moreover, when trajectories are close to $\bar \gamma$, $\|
(s,\dot s)\|$ is small, implying that $(q,\dot q)$ is close to the
original constraint manifold $\Gamma$, as argued above.

The two control specifications we have defined so far, namely the
asymptotic stabilization of both $\bar \Gamma$ and  $\bar \gamma$, are
somewhat related to one another in that $\bar \gamma \subset \bar
\Gamma$.  Indeed, on $\bar \gamma$ one has that $(q,\dot q)
=(\sigma(\theta),\sigma'(\theta)\dot \theta) \in \Gamma$ and $(s,\dot
s)=(0,0)$, which readily implies that $(q,\dot q,s,\dot s) \in \bar
\Gamma$.

In summary, we have formulated the following

\textbf{\vhc-based orbital stabilization
  problem.} \label{prob:stab:VHCorbitStab} Find a smooth control law
for system~\eqref{eq:stab:ELsys} augmented with the double-integrator
$\ddot s =v$ 
that asymptotically stabilizes both  sets $\bar \gamma \subset
\bar \Gamma$.

We recall from the foregoing discussion that
the asymptotic stabilization of $\bar \Gamma$ corresponds to the
enforcement of the perturbed \vhc $h^s(q)=0$. Since $(s,\dot s)=(0,0)$
on $\bar \gamma$, near $\bar \gamma$ the Hausdorff
distance\footnote{The Hausdorff distance between two sets measures how
  far the two sets are from each other.} between the set $\bar \Gamma$
and the original constraint manifold $\Gamma\times \{(s,\dot
s)=(0,0)\}$ is small. Considering the fact that $h(q)=0$ embodies a
useful constraint that we wish to hold during the transient, the
philosophy of the \vhc-based orbital stabilization problem is to
preserve as much as possible the beneficial properties of the original
\vhc $h(q)=0$, while simultaneously stabilizing the closed orbit
$\gamma$ corresponding to a desired repetitive motion.

\textbf{Solution steps.} Our solution to the \vhc-based orbital
stabilization problem unfolds in three  steps:
\begin{enumerate}

\item We present a technique to parameterize the \vhc $h(q)=0$ with the
  output of a double integrator, giving rise to a dynamic \vhc
  $h^s(q)=0$ with associated constraint manifold $\bar \Gamma$. We
  show that if the original \vhc is regular, so too is its dynamic
  counterpart for small values of the double integrator output
  (Proposition~\ref{prop:stab:augmented}). Moreover, if the original
  constraint manifold $\Gamma$ is stabilizable, so too is the
  perturbed manifold $\bar \Gamma$
  (Proposition~\ref{prop:stab:stabilizability}). We derive the reduced
  dynamics on this manifold, which are now affected by the input $v$
  of the double integrator.

\item We develop a general result for control-affine systems
  (Theorem~\ref{thm:transverse_linearization}) relating the
  exponential stabilizability of a closed orbit to the controllability
  of a linear periodic system, for which we give an explicit
  representation. Leveraging this result, we design the input of the
  double-integrator, $v$, to exponentially stabilize the orbit
  relative to the manifold $\bar \Gamma$.

\item We put together the controller enforcing the dynamic \vhc in
  Step 1 with the controller stabilizing the orbit in Step 2 and show
  that the resulting controller solves the \vhc-based orbital
  stabilization problem (Theorem~\ref{thm:main_result}).
\end{enumerate}
%
\section{Step 1: Making the \vhc dynamic}\label{sec:stab:dynConstr}  
In this section we present the notion of dynamic \vhcs. We begin by
augmenting the dynamics in~\eqref{eq:stab:ELsys} with a
double-integrator, to obtain the augmented system
\begin{equation}\label{eq:stab:sys_augmented}
\begin{aligned}
D(q) \ddot q +  C(q,\dot q) \dot q + \nabla P(q) &= B(q) \tau, \\
 \ddot s &= v.
\end{aligned}
\end{equation}
Henceforth, we use overbars to distinguish objects associated with the
augmented control system \eqref{eq:stab:sys_augmented} from those
associated with~\eqref{eq:stab:ELsys}. Accordingly, we define
$\bar{q}:=(q,s)$, $\dot{\bar{q}}:=(\dot{q},\dot{s})$,
$\bar{\cQ}:=\big\{(q,s): q\in\cQ,\, s\in \mathcal{I}\big\}$.
\begin{definition}\label{defn:regular:dynamic}
Let $h(q)=0$ be a regular \vhc of order $n-1$ for
system~\eqref{eq:stab:ELsys}.  A \textbf{dynamic \vhc based on
  $h(q)=0$} is a relation $h^s(q)=0$ such that the map $(s,q)\mapsto
h^s(q)$ is $C^2$, $h^0(q)=h(q)$, and the parameter $s$ satisfies the
differential equation $\ddot{s} = v$ in~\eqref{eq:stab:sys_augmented}.

The dynamic \vhc $h^s(q)$ is \textbf{regular}
for~\eqref{eq:stab:sys_augmented} if there exists an open interval
$\mathcal{I}\subset \Re$ containing $s=0$ such that, for all $s\in
\mathcal{I}$ and all $v\in \Re$, system~\eqref{eq:stab:sys_augmented}
with input $\tau$ and output $e=h^{s}(q)$ has vector relative degree
$\{2,\cdots,2\}$.

The dynamic \vhc $h^s(q)=0$ is \textbf{stabilizable}
for~\eqref{eq:stab:sys_augmented} if there exists a smooth feedback
$\tau(q,\dot{q},s,\dot{s},v)$ such that the manifold
\begin{equation}
\label{eq:stab:Gammabar0}
\bar{\Gamma} :=  \{(q,\dot{q},s,\dot{s}): 
h^s(q)=0,\, \partial_q h^s\dot{q} + \partial_s h^s \dot{s}=0\},
\end{equation}
is asymptotically stable for the closed-loop system.
\end{definition}

The reason for parameterizing the \vhc with the output of a double
integrator is to guarantee that the input $v$ of the double integrator
appears after taking two derivatives of the output function
$e=h^{s}(q)$. The regularity property of $h^{s}(q)$ in the foregoing
definition means that, upon calculating the second derivative of
$e=h^s(q)$ along the vector field in~\eqref{eq:stab:sys_augmented},
the control input $\tau$ appears nonsingularly, i.e.,
\[
\ddot{e} = (\star) + A^{s}(q)\tau + B^{s}(q)v, 
\]
where $A^s$ and $B^s$ are suitable matrices, and $A^s$ is invertible
for all $q\in (h^{s})^{-1}(0)$ and all $s\in \cI$.

\begin{figure}[ht]
\psfrag{h}[c]{$h(q)=0$}
\psfrag{j}{$h^s(q)=0$}
\psfrag{L}{$Ls$}
\psfrag{x}{$q_1$}
\psfrag{y}{$q_n$}
\centerline{\includegraphics[width=.4\textwidth]{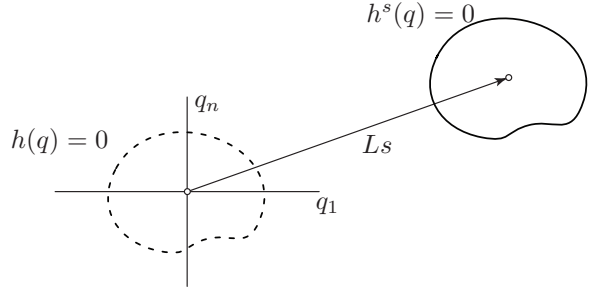}}
\caption{Geometric interpretation of the dynamic
  \vhc~\eqref{eq:stab:dynVHC}. The dashed closed curve represents
    the original \vhc $h(q)=0$, while the solid curve represents the
    dynamic \vhc $h^s(q)=0$. The new configuration variable $s$
      parametrizing the \vhc has the effect of translating the
      original \vhc curve in the direction spanned by the vector $L$.}
\label{fig:stab:dynVHC}
\end{figure}

Given a regular \vhc~$h(q)=0$, a possible way to generate a dynamic \vhc
based on $h(q)=0$ is to translate the curve $h^{-1}(0)$ by an
amount proportional to $s \in \Re$. Other choices are of course possible, but
this one has the benefit of allowing for simple expressions in the
derivations that follow.  We thus consider the following
one-parameter family of mappings
\begin{equation} \label{eq:stab:dynVHC}
h^s(q):= h(q-Ls), 
\end{equation}
where $L \in \Re^n$ is a non-zero constant vector.  The zero level set
of each family member in~\eqref{eq:stab:dynVHC} is
$(h^s)^{-1}(0)=\{q+Ls: q\in h^{-1}(0)\}$, a
translation\footnote{Recall that $q$ is a $n$-tuple whose $i$-th
  element, $q_i$, is either a real number or an element of
  $[\Re]_{T_i}$. In the latter case, the sum $q_i + L_i s$ is to be
  understood as sum modulo $T_i$.}
 of $h^{-1}(0)$
by the vector $Ls$ (see Figure~\ref{fig:stab:dynVHC}). If
$\sigma:[\Re]_{\Tconstraint}\rightarrow \cQ$ is a regular parameterization of the
curve $h^{-1}(0)$, a regular parameterization of the zero level
set of each family member in~\eqref{eq:stab:dynVHC} is
$\sigma^{s}(\theta)= \sigma(\theta)+Ls$.  In an analogous manner, the
constraint manifold $\bar{\Gamma}$ in~\eqref{eq:stab:Gammabar0} is the
translation of $\Gamma$ by the vector $\col(Ls,\, L\dot{s})$,
\begin{equation}\label{eq:stab:Gammabar}
\begin{aligned}
\bar{\Gamma} &= \{ (q,\dot{q},s,\dot{s}): h(q-Ls)=0,\,
dh_{q-Ls}(\dot{q}-L\dot{s})=0\} \\
&= \{(q,\dot{q},s,\dot{s}): (q-Ls,\dot{q}-L\dot{s})\in \Gamma\}.
\end{aligned}
\end{equation}
In the  augmented coordinates, the closed orbit we wish to stabilize is 
\begin{equation}\label{eq:gammabar}
\bar \gamma = \{(q,\dot q,s,\dot s) : s=\dot s=0, \, q =
\sigma(\theta), \, \dot q = \sigma'(\theta) \dot \theta, \,
(\theta,\dot \theta) \in [\Re]_{\Tconstraint} \times \Re, \,
E(\theta,\dot \theta)=E_0\}.
\end{equation}
The next two propositions show that if $h(q)=0$ is regular and
stabilizable, so too is its dynamic counterpart $h(q-Ls)=0$. Their
proofs are in Appendix~\ref{sec:appendixA}.

\begin{prop}\label{prop:stab:augmented}
	If $h(q)=0$ is a regular \vhc of order $n-1$
        for~\eqref{eq:stab:ELsys}, then for any $L \in \Re^n$ the
        dynamic \vhc $h(q-Ls)=0$ is regular for the augmented
        system~\eqref{eq:stab:sys_augmented}.
\end{prop}


\begin{prop}\label{prop:stab:stabilizability}
If $h(q)=0$ is a regular \vhc of order $n-1$ for~\eqref{eq:stab:ELsys}
satisfying the stabilizability condition~\eqref{eq:classK}, then for
any $L\in \Re^n$, the dynamic \vhc $h(q-Ls)=0$ is stabilizable in the
sense of Definition~\ref{defn:regular:dynamic}, and a feedback
stabilizing the constraint manifold $\bar \Gamma$
in~\eqref{eq:stab:Gammabar} is $\tau= \tau^\star(q,\dot q,s,\dot s,v)$
given by
\begin{equation} \label{eq:stab:tau_s_method}
\tau^\star(q,\dot q,s,\dot s,v)=\big(A^s(q)\big)^{-1}\Big\{dh_{q-Ls}\,
D^{-1}(q)\big[C(q,\dot{q})\dot{q}+\nabla P(q)\big]+ dh_{q-Ls}\,
Lv-\mathcal{H}(q,\dot{q},s,\dot{s})-k_p e-k_d \dot{e}\Big\},
\end{equation}
where $e=h(q-Ls)$, $\dot e = dh_{q-Ls} \, (\dot q - L \dot s)$,
$A^{s}(q)=dh_{q-Ls}\, D^{-1}(q)B(q)$, $k_p$, $k_d>0$, and
$\cH=\col(\cH_1,\ldots,\cH_{n-1})$, $\cH_i=
(\dot{q}-L\dot{s}){\trans} \Hess(h_i)|_{q-Ls}
  (\dot{q}-L\dot{s})$.
\end{prop}

%
%
%
Next, we find the reduced dynamics of the augmented
system~\eqref{eq:stab:sys_augmented} with
feedback~\eqref{eq:stab:tau_s_method} on the manifold $\bar \Gamma$
in~\eqref{eq:stab:Gammabar}. To this end, we
left-multiply~\eqref{eq:stab:sys_augmented} by the left annihilator
$B^\perp$ of $B$ and evaluate the resulting equation on $\bar{\Gamma}$
by setting
\[
q=\sigma(\theta)+Ls, \ \dot q = \sigma'(\theta)\dot \theta+L \dot s,
\ \ddot q = \sigma'(\theta)\ddot \theta + \sigma''(\theta) \dot
\theta^2 + L v.
\]
By so doing, one obtains:
\begin{equation}
\label{eq:stab:reducedDyn_extended}
\begin{aligned}
& \ddot{\theta} = \Psi_1^s(\theta)+\Psi_2^s(\theta)\dot{\theta}^2+\Psi_3^s(\theta)\dot{\theta}\dot{s}+
\Psi_4^s(\theta)\dot{s}^2+\Psi_5^s(\theta)v, \\
&
\ddot{s} = v,
\end{aligned}
\end{equation}
where 
\begin{equation}\begin{aligned}
\Psi_1^s(\theta) &= -\frac{B^{\perp} \nabla P}{B^{\perp}D\sigma^{\prime}}\bigg|_{q=\sigma(\theta)+Ls},
\\
\Psi_2^s(\theta)  &= -\frac{B^{\perp}D\sigma^{\pprime} +\sum\limits_{i=1}^{n} B_i^{\perp}\sigma^{\prime{\trans}}Q_i\sigma^{\prime}}{B^{\perp}D\sigma^{\prime}}\Bigg|_{q=\sigma(\theta)+Ls},\\
\Psi_3^s(\theta) &= -\frac{2\sum\limits_{i=1}^{n} B_i^{\perp}\sigma^{\prime{\trans}}Q_i L}{B^{\perp}D\sigma^{\prime}}\Bigg|_{q=\sigma(\theta)+Ls},\\
\Psi_4^s(\theta) &= -\frac{\sum\limits_{i=1}^{n} B_i^{\perp}L{\trans}Q_i L}{B^{\perp}D\sigma^{\prime}}\Bigg|_{q=\sigma(\theta)+Ls},
\\
\Psi_5^s(\theta) &= -\frac{B^{\perp} D L}{B^{\perp}D\sigma^{\prime}}\bigg|_{q=\sigma(\theta)+Ls}.
  \end{aligned}
\end{equation}
The two second-order differential
equations~\eqref{eq:stab:reducedDyn_extended} will be henceforth
referred to as the {\bf extended reduced dynamics} induced by the
dynamic \vhc $h(q-Ls)=0$. They represent the motion of the mechanical
system~\eqref{eq:stab:ELsys} on the constraint manifold $\bar \Gamma$
in~\eqref{eq:stab:Gammabar}.  Their restriction to the plane $\{s=\dot
s=0\}$ coincides with the reduced
dynamics~\eqref{eq:constrained_dynamics}. The state $(\theta,\dot
\theta,s,\dot s) \in [\Re]_{\Tconstraint} \times \Re \times \Re \times
\Re$ represents global coordinates for $\bar \Gamma$. Using these
coordinates, and with a slight abuse of notation, the closed orbit
$\bar \gamma$ in~\eqref{eq:gammabar} is given by
\begin{equation}\label{eq:gammabar2}
\bar \gamma =\{ (\theta,\dot \theta,s,\dot s) \in [\Re]_{\Tconstraint} \times \Re
\times \Re\times \Re: E(\theta,\dot \theta) = E_0, \, s=\dot s=0\}.
\end{equation}
This is the set we will stabilize next.

\section{Step 2: Linearization along the closed orbit}\label{sec:stab:s_control}

The objective now is to design the control input $v$ in the extended
reduced dynamics~\eqref{eq:stab:reducedDyn_extended} so as to
stabilize the closed orbit $\bar \gamma$ in~\eqref{eq:gammabar2}. We
will do so by adopting the philosophy of Hauser et
al. in~\cite{HauChu94} that relies on an implicit representation of
the closed orbit to derive the so-called transverse linearization
along $\bar \gamma$. Roughly speaking, this is the linearization along
$\bar \gamma$ of the components of the dynamics that are transversal
to $\bar \gamma$. Hauser's approach generalizes classical results of
Hale~\cite[Chapter VI]{Hale1980},  requiring a moving orthonormal
frame. The insight in~\cite{HauChu94} is that orthogonality is not
needed, transversality is enough. This insight allowed Hauser et
al. in~\cite{HauChu94} to derive a normal form analogous to that
in~\cite[Chapter VI]{Hale1980}, but calculated directly from an
implicit representation of the orbit. We shall use the same idea in
the theorem below.

We begin by enhancing the results of~\cite{Hale1980,HauChu94} in two
directions. First, while~\cite{Hale1980,HauChu94} require the
knowledge of a periodic solution, we only require a parameterization
of $\bar \gamma$ (something that is readily available in the setting
of this paper, while the solution is not). Second,
while~\cite{Hale1980,HauChu94} deals with dynamics without inputs, we
provide a necessary and sufficient criterion for the exponential
stabilizability of the orbit.

{\bf A general result.} Our first result is a necessary and sufficient
condition for a closed orbit to be exponentially stabilizable. This
result is of considerable practical use, and is of independent
interest. 

Consider a control-affine system
\begin{equation}\label{eq:control_affine}
\dot x = f(x) +g(x) u,
\end{equation}
with state $x \in \cX$, where $\cX$ is a closed embedded submanifold
of $\Re^n$, and control input $u \in \Re^m$. A closed orbit $\gamma$
is {\bf exponentially stabilizable} for~\eqref{eq:control_affine} if
there exists a locally Lipschitz continuous feedback $u^\star(x)$ such
that the set $\gamma$ is exponentially stable for the closed-loop
system $\dot x = f(x) +g(x) u^\star(x)$, i.e., there exist
$\delta,\lambda,M>0$ such that for all $x_0 \in \cX$ such that
$\|x_0\|_\gamma < \delta$, the solution $x(t)$ of the closed-loop
system satisfies $\|x(t)\|_\gamma \leq M \|x_0\|_\gamma e^{-\lambda
  t}$ for all $t \geq 0$.  Note that if $\gamma$ is exponentially
stable, then $\gamma$ is asymptotically stable.

Let $T$ be a positive real number. A linear $T$-periodic system
$dx/dt= A(t)x$, where $A(\cdot)$ is a continuous and $T$-periodic
matrix-valued function, is {\bf asymptotically stable} if all its
characteristic multipliers lie in the open unit disk. A linear
$T$-periodic control system
\begin{equation}\label{eq:LTP}
\frac{dx}{dt} = A(t) x + B(t) u,
\end{equation}
where $A(\cdot)$ and $B(\cdot)$ are continuous and $T$-periodic
matrix-valued functions, is {\bf stabilizable} (or the pair
$(A(\cdot),B(\cdot))$ is stabilizable) if there exists a continuous
and $T$-periodic matrix-valued function $K(\cdot)$ such that $\dot x =
(A(t) + B(t) K(t))x$ is asymptotically stable. In this case, we say
that the feedback $u=K(t) x$ {\bf stabilizes}
system~\eqref{eq:LTP}. The notion of stabilizability can be
characterized in terms of the characteristic multipliers of $A(\cdot)$
(see, e.g.~\cite{bittanti1991periodic}).

\begin{thm}\label{thm:transverse_linearization}
Consider system~\eqref{eq:control_affine}, where $f$ is a $C^1$ vector
field and $g$ is locally Lipschitz continuous on $\cX$. Let $\gamma
\subset \cX$ be a closed orbit of the open-loop system $\dot x =
f(x)$, and let $\s \mapsto \varphi(\s)$, $[\Re]_T \to \cX$, be a
regular parameterization of $\gamma$.  Finally, let $H:\cX \to
\Re^{n-1}$ be an implicit representation of $\gamma$ with the
properties that $H$ is $C^1$, $\rank dH_x =n-1$ for all $x \in
H^{-1}(0)$, and $H^{-1}(0)=\gamma$.

\begin{enumerate}[(a)]

\item The orbit $\gamma$ is exponentially stabilizable
  for~\eqref{eq:control_affine} if, and only if, the linear
  $T$-periodic control system on $\Re^{n-1}$
\begin{equation}\label{eq:transverse_linearization:general}
\begin{aligned}
&\dot \z = A(t) z + B(t) u \\
&A(t)=\frac{\| \varphi'(t)\|^2 }{\langle f(\varphi(t)), \varphi'(t)
    \rangle} \left[ (dL_f H)_{\varphi(t)} dH^\dagger_{\varphi(t)}
    \right] \\
&B(t)=\frac{\| \varphi'(t)\|^2 }{\langle f(\varphi(t)), \varphi'(t)
    \rangle} \Big[ L_g H(\varphi(t)) \Big],
\end{aligned}
\end{equation}
is stabilizable.

\item If a $T$-periodic feedback $u=K(t) \z$, with $K(\cdot)$ continuous and
  $T$-periodic, stabilizes the $T$-periodic
  system~\eqref{eq:transverse_linearization:general}, then for any
  smooth map $\pi : \cU \to [\Re]_T$, with $\cU$ a neighbourhood of
  $\gamma$ in $\cX$ and $\pi$ such that $\pi|_\Gamma= \varphi^{-1}$,
  the feedback
\begin{equation}\label{eq:ustar}
u^\star(x)=K(\pi(x)) H(x)
\end{equation}
exponentially stabilizes the closed orbit $\gamma$
for~\eqref{eq:control_affine}.
\end{enumerate}
\end{thm}
The proof of Theorem~\ref{thm:transverse_linearization} is found in
Appendix~\ref{sec:appendixB}.
%
%
%
\begin{rem}
Concerning the existence of the function $H$ in the theorem statement,
since closed orbits of smooth dynamical systems are diffeomorphic to
the unit circle $\Se^1$, it is always possible to find a function $H$
satisfying the assumptions of the theorem.  This well-known fact is
shown, e.g., in~\cite[Proposition 1.2]{HauChu94}. Since the conditions
of the theorem are necessary and sufficient, the result is independent
of the choice of $H$.  As for the existence of the function $\pi : \cU
\to [\Re]_T$, this function can be constructed by picking $\cU$ to be
a tubular neighbourhood of $\gamma$.  Then there exists a smooth
retraction $r: \cU \to \gamma$. The function $\pi = \varphi^{-1}\circ
r$ has the desired properties. If $\varphi(\cdot)$ is a
$T$-periodic solution of $\dot x=f(x)$, rather than just a regular
parametrization of $\gamma$, then we have $\varphi'(t) =
f(\varphi(t))$, and the scalar coefficient
in~\eqref{eq:transverse_linearization:general}, $\| \varphi'(t)\|^2
/ \langle f(\varphi(t)),\varphi'(t) \rangle$, is identically equal
to one.  \hfill $\triangle$
\end{rem}
\begin{rem}
Theorem~\ref{thm:transverse_linearization} establishes the equivalence
between the exponential stabilizability of the closed orbit $\gamma$
and the stabilizability of the linear periodic
system~\eqref{eq:transverse_linearization:general}, the so-called
transverse linearization.  The equivalence between these two concepts
is not new, it is essentially contained in the results
of~\cite[Chapter VI]{Hale1980} and~\cite{HauChu94}. What is new in
Theorem~\ref{thm:transverse_linearization}, and of considerable
practical interest, is the fact that it provides an explicit
expression for the transverse linearization that can be computed using
{\em any} regular parametrization $\varphi$ of $\gamma$ and {\em any}
implicit representation $H$ of $\gamma$ whose Jacobian matrix has full
rank on $\gamma$. In contrast to the above, the methods
in~\cite{Hale1980} and~\cite{HauChu94} rely on the knowledge of a
periodic open-loop solution of~\eqref{eq:control_affine} generating
$\gamma$ and do not give an explicit expression for
$(A(\cdot),B(\cdot))$. We also mention that Hauser's notion of
transverse linearization was applied in~\cite{ShiPerWit05} to a
special class of Euler-Lagrange systems, once again requiring the knowledge of a
periodic solution. Moreover, in~\cite{ShiFreGus10}, the authors do
give an explicit expressions for the transverse linearization
$(A(\cdot),B(\cdot))$, but one that is only applicable to a class of
Euler-Lagrange systems, while the expressions in
Theorem~\ref{thm:transverse_linearization} are applicable to
arbitrary vector fields.  \hfill $\triangle$
\end{rem}

\begin{rem} In the special case of systems without control (i.e.,
  $g(x)=0$ in~\eqref{eq:control_affine}),
  Theorem~\ref{thm:transverse_linearization} implies that the closed
  orbit $\gamma$ is exponentially stable if and only if the origin of
  the linear periodic system $\dot z = A(t) z$, with $A(t)$ given
  in~\eqref{eq:transverse_linearization:general}, is asymptotically
  stable or, equivalently,
  multipliers of $A(t)$ have magnitude $<1$.  This result is to be
  compared to the Poincar\'e stability theorem also known as the
  Andronov-Vitt theorem (AVT)~\cite{AndVit33} (see also~\cite[Chapter
    VI, Theorem~2.1]{Hale1980}), stating that if $\varphi(t)$ is a
  $T$-periodic solution of the dynamical system $\dot x=f(x)$, then
  the closed orbit $\gamma =\image(\varphi)$ is orbitally stable if
  the characteristic multipliers of the {\em variational equation}
  $\dot x = (df_{\varphi(t)}) x$ are $\{1, \mu_1,\ldots, \mu_{n-1}\}$,
  with $|\mu_i|<1$, $i=1,\ldots, n-1$.  The link between the
  Andronov-Vitt theorem and Theorem~\ref{thm:transverse_linearization}
  is that the complex numbers $\{\mu_1,\ldots,\mu_{n-1}\}$ in the AVT
  are the characteristic multipliers of the $n-1 \times n-1$ matrix
  $A(t)$ in Theorem~\ref{thm:transverse_linearization}.  There are,
  however, two important differences.  First, as already pointed out,
  Theorem~\ref{thm:transverse_linearization} provides a means to
  directly calculate $\{\mu_1,\ldots,\mu_{n-1}\}$ without the need to
  know a periodic solution of the differential equation $\dot x
  =f(x)$. Rather, any parametrization $\varphi(\s)$ of $\gamma$
  suffices. Furthermore, even if a periodic solution is available, the
  transverse linearization $A(t)$
  in~\eqref{eq:transverse_linearization:general} differs from the one
  in the literature (see equation (1.9) in~\cite[Chapter
    VI]{Hale1980}) because it does not rely an {\em orthonormal}
  moving frame\footnote{For the method of orthonormal moving frames,
    the reader may also consult~\cite{Ura58}.} (the matrix $Z(\s)$
  in~\cite{Hale1980}). Indeed, the columns of the differential of $H$
  appearing in the definition of $A(t)$ span the plane orthogonal to
  the tangent vector to the orbit, $\varphi'(t)$, but they do not
  necessarily form an orthonormal frame.

  Finally, in the context of stability of orbits of dynamical systems,
  we mention the work of Demidovich in~\cite{Dem68} which generalized
  the work of Andronov-Vitt for the stability of not necessarily closed
  orbits, and the work in~\cite{Leo06} which further generalized
  Demidovich's work. When specialized to closed orbits, the results
  in~\cite{Dem68,Leo06} differ from
  Theorem~\ref{thm:transverse_linearization} in the same way that the
  AVT does.  \hfill $\triangle$
\end{rem}

{\bf Design of the $T$-periodic feedback matrix $K$.} Once it is
established that the pair $(A(\cdot),B(\cdot))$
in~\eqref{eq:transverse_linearization:general} is stabilizable, the
design of the $T$-periodic feedback matrix $K(\cdot)$ in part~(b) of
the theorem can be carried out by solving the periodic Riccati
equation for a $T$-periodic $\Pi: \Re \to \Re^{n-1 \times n-1}$:
\begin{equation} \label{eq:PRE}
\begin{aligned}
-\frac{d \Pi}{dt}=A(t){\trans}\Pi(t)+\Pi(t)A(t)
-\Pi(t)B(t)R^{-1}B(t){\trans}\Pi(t)+Q(t).
\end{aligned}
\end{equation}
where $R(\cdot)=R(\cdot)\trans$ is a positive definite continuous
$T$-periodic matrix-valued function and $Q(\cdot) = Q(\cdot)\trans$ is
a positive definite continuous $T$-periodic matrix-valued function,
and setting
\begin{equation}\label{eq:K}
K(t) = -\frac 1 R B(t)\trans \Pi(t).
\end{equation}
Theorem~6.5 in~\cite{bittanti1991periodic} states that if, and only
if, $(A(\cdot),B(\cdot))$ is stabilizable and
$(Q^{1/2}(\cdot),A(\cdot))$ is detectable (this latter condition is
satisfied, e.g., by letting $Q$ be the identity matrix) then the
Riccati equation~\eqref{eq:PRE} has a unique positive semidefinite
$T$-periodic solution $\Pi(\cdot)$ and the feedback $u = K(t) z$, with
$K(\cdot)$ given in~\eqref{eq:K}, stabilizes
system~\eqref{eq:transverse_linearization:general}. Once this is done,
the feedback $u^\star(x)$ in~\eqref{eq:ustar} exponentially stabilizes
the closed orbit $\gamma$.

{\bf Application to extended reduced dynamics.}
We now apply Theorem~\ref{thm:transverse_linearization} to the extended reduced
dynamics~\eqref{eq:stab:reducedDyn_extended} with the objective of
stabilizing the closed orbit $\bar \gamma$
in~\eqref{eq:gammabar2}. Here we have $x=(\theta,\dot \theta,s,\dot
s)$, and the implicit representation of $\bar \gamma$
\[
H(x) = ( E(\theta,\dot \theta)-E_0, s,\dot s).
\]
Leveraging the parameterizations of closed orbits presented in
Section~\ref{sec:stab:prelim}, the parameterization of $\bar \gamma$
in $(\theta,\dot \theta,s,\dot s)$-coordinates has the form
$[\Re]_{\Tgamma} \to [\Re]_{\Tconstraint} \times \Re \times \Re \times
\Re$, $\s \mapsto (\varphi_1(\s), \varphi_2(\s),0,0)$, with
$\varphi_1$, $\varphi_2$ given by~\eqref{eq:param_rotation} and
$\Tgamma = \Tconstraint$ if $\gamma$ is a rotation, and
by~\eqref{eq:param_oscillation} and $\Tgamma = 2\pi$ if $\gamma$ is an
oscillation. Applying Theorem~\ref{thm:transverse_linearization} to
system~\eqref{eq:stab:reducedDyn_extended}, we get the following
$\Tgamma$-periodic linear system
\begin{equation}\label{eq:transverse_linearization}
\dot z = 
\begin{bmatrix}
0 & a_{12}(t) & a_{13}(t)\\
0 & 0 & 1 \\
0 & 0 & 0
\end{bmatrix} z + 
\begin{bmatrix}
b_1(t) \\
0 \\
1
\end{bmatrix}
v,
\end{equation}
where 
\begin{equation}\label{eq:A_and_b_components}
\begin{aligned}
& a_{12}(t) = \eta(t) M(\varphi_1(t)) \varphi_2(t) \big[\partial_{\z_2}
    \Psi_1^{\z_2}(\varphi_1(t)) + \partial_{\z_2}
    \Psi_2^{\z_2}(\varphi_1(t)) \varphi_2^2(t)\big]\big|_{\z_2=0}, \\
& a_{13}(t) = \eta(t) M(\varphi_1(t)) \varphi_2^2(t) \Psi_3^0(\varphi_1(t)), \\
& b_1(t) = \eta(t) M(\varphi_1(t)) \varphi_2(t) \Psi_5^0(\varphi_1(t)), \\
& \eta(t) =
  \frac{(\varphi_1'(t))^2+(\varphi_2'(t))^2}{\varphi_1'(t)\varphi_2(t) +
    \varphi_2'(t) [\Psi_1(\varphi_1(t)) + \Psi_2(\varphi_1(t))
      \varphi_2^2(t) ]}.
\end{aligned}
\end{equation}
Assuming that system~\eqref{eq:transverse_linearization} is
stabilizable, then we may find the unique positive semidefinite
solution of the periodic Riccati equation~\eqref{eq:PRE} to get the
matrix-valued function $K(\cdot)$
in~\eqref{eq:K}. Theorem~\ref{thm:transverse_linearization} guarantees
that the controller
\begin{equation}\label{eq:vstar}
v=\bar v(\theta,\dot \theta,s,\dot s) = K(\pi(\theta,\dot \theta,s,\dot
s)) \begin{bmatrix} E(\theta,\dot \theta) - E_0 \\ s \\ \dot s
\end{bmatrix}
\end{equation}
exponentially stabilizes the orbit $\bar \gamma$
in~\eqref{eq:gammabar2} for the extended reduced
dynamics~\eqref{eq:stab:reducedDyn_extended}.

It remains to find an explicit expression for the map $\pi$. If
$\gamma$ is a rotation, then in light of the
parameterization~\eqref{eq:param_rotation}, we may set
\[
\pi(\theta,\dot \theta,s,\dot s) = \theta.
\]
Else, if $\gamma$ is an oscillation,
using~\eqref{eq:param_oscillation} we set
\[
\pi(\theta,\dot \theta,s,\dot s) = \text{atan2} (T(\theta) \dot \theta,\theta-C),
\]
where $\text{atan2}(\cdot,\cdot)$ is the four-quadrant arctangent
function such that $\text{atan2}(\sin(\alpha),\cos(\alpha)) = \alpha$
for all $\alpha \in (-\pi,\pi)$.

\section{Step 3: Solution of the \vhc-based orbital stabilization problem}\label{sec:stab:sol}

In Section~\ref{sec:stab:dynConstr}, we designed the feedback
$\tau^\star$ in~\eqref{eq:stab:tau_s_method} to asymptotically
stabilize the constraint manifold $\bar \Gamma$ associated with the
dynamic \vhc $h(q-L s)=0$. In Section~\ref{sec:stab:s_control}, we
designed the feedback $\bar v$ in~\eqref{eq:vstar} for the double
integrator $\ddot s =v$ rendering the closed orbit $\bar \gamma$
exponentially stable relative to $\bar \Gamma$ (i.e., when initial
conditions are on $\bar \Gamma$).  There are two things left to do in
order to solve the \vhc-based orbital stabilization problem. First, in
order to implement the feedback $\bar v$ in~\eqref{eq:vstar}, we need
to relate the variables $(\theta,\dot \theta)$ to the state $(q,\dot
q)$. Second, we need to show that the asymptotic stability of $\bar
\Gamma$ and the asymptotic stability of $\bar \gamma$ relative to
$\bar \Gamma$ imply that $\bar \gamma$ is asymptotically stable.

To address the first issue, we leverage the fact that, since
$h^{-1}(0)$ is a closed embedded submanifold of $\cQ$,
by~\cite[Proposition~6.25]{Lee13} there exists a neighbourhood $\cW$ of
$h^{-1}(0)$ in $\cQ$ and a smooth retraction of $\cW$ onto
$h^{-1}(0)$, i.e., a smooth map $r : \cW \to h^{-1}(0)$ such that
$r|_{h^{-1}(0)}$ is the identity on $h^{-1}(0)$. Define $\Theta:\cW
\to [\Re]_T$ as $\Theta = \sigma^{-1} \circ r$. By construction,
$\Theta|_{h^{-1}(0)} =\sigma^{-1}$. In other words, for all $q \in
h^{-1}(0)$, $\Theta(q)$ gives that unique value of $\theta \in
[\Re]_T$ such that $q = \sigma(\theta)$. Using the function $\Theta$,
we now define an extension of $\bar v$ from $\bar \Gamma$ to a
neighborhood of $\bar \Gamma$ as follows
\begin{equation}\label{eq:vstar:extended}
v^\star(q,\dot q,s,\dot s) = \bar v(\theta,\dot \theta,s,\dot s)
\big|_{(\theta,\dot \theta) = (\Theta(q),d\Theta_q \dot q)}.
\end{equation}
We are now ready to solve the \vhc-based orbital stabilization
problem.
\begin{thm}\label{thm:main_result}
Consider system~\eqref{eq:stab:ELsys} and let $h(q)=0$ be a regular
\vhc of order $n-1$. Let $\sigma:[\Re]_\Tconstraint \to \cQ$ be a
regular parametrization of $h^{-1}(0)$ and consider the following
assumptions:

\begin{enumerate}[(a)]

\item The \vhc $h(q)=0$ satisfies the stabilizability
  condition~\eqref{eq:classK}.

\item The \vhc $h(q)=0$ induces Lagrangian reduced dynamics as per
  Proposition~\ref{prop:ESAIM}.

\item For a closed orbit $\gamma$ of the reduced dynamics given in
  implicit form as $\gamma = \{ (\theta,\dot \theta) \in
  [\Re]_\Tconstraint \times \Re: E(\theta,\dot \theta) =E_0\}$,
  consider one of the regular parametrizations $[\Re]_{\Tgamma} \mapsto
  [\Re]_\Tconstraint \times \Re\times \Re \times \Re$ discussed in
  Section~\ref{sec:stab:s_control}. Assume that the $\Tgamma$-periodic
  system~\eqref{eq:transverse_linearization}-\eqref{eq:A_and_b_components}
  is stabilizable.

\end{enumerate}

Under the assumptions above, let $Q(\cdot)=Q(\cdot)\trans$ be a
positive definite $\Tgamma$-periodic $\Re^{3\times 3}$-valued function
such that $(Q^{1/2},A)$ is detectable, and pick any $R>0$.  The smooth
  dynamic feedback
\[
\begin{aligned}
& \tau = \tau^\star(q,\dot q,s,\dot s,v^\star(q,\dot q,s,\dot s))\\
& \ddot s = v^\star(q,\dot q,s,\dot s),
\end{aligned}
\]
with $\tau^\star$ defined in~\eqref{eq:stab:tau_s_method}, $v^\star$
defined in~\eqref{eq:vstar},~\eqref{eq:vstar:extended}, and where
$K(\cdot)$ in~\eqref{eq:K} results from the solution of the
$\Tgamma$-periodic Riccati equation~\eqref{eq:PRE}, asymptotically
stabilizes both sets $\bar \gamma \subset \bar \Gamma$ given
in~\eqref{eq:stab:Gammabar},~\eqref{eq:gammabar}.
\end{thm}
\begin{pf} 
Propositions~\ref{prop:stab:augmented}
and~\ref{prop:stab:stabilizability} establish that the
feedback~\eqref{eq:stab:tau_s_method} stabilizes the set $\bar
\Gamma$. Theorem~6.5 in~\cite{bittanti1991periodic} establishes that
$K(\cdot)$ in~\eqref{eq:K} is well-defined, and
Theorem~\ref{thm:transverse_linearization} establishes that $v^\star$
in~\eqref{eq:vstar:extended} stabilizes $\bar \gamma$ relative to
$\bar \Gamma$. Since $\bar \gamma$ is a compact set, the reduction
theorem for stability of compact sets
in~\cite{seibert-1995},~\cite{el-hawwary-2013} implies that $\bar
\gamma$ is asymptotically stable for the closed-loop system.
\hfill\(\qed\) \end{pf}

\begin{figure}[ht]
\psfrag{M}[cb]{Mechanical system} \psfrag{V}[cb]{\vhc stabilizer}
\psfrag{O}[cb]{Orbit stabilizer} \psfrag{s}{$s$} \psfrag{S}{$\dot s$}
\psfrag{v}{$v$} \psfrag{t}[l]{$(\theta,\dot \theta)$}
\psfrag{q}{$(q,\dot q)$} \psfrag{T}{$\tau$}
\psfrag{I}[c]{$\bigintsss$}
\centerline{\includegraphics[width=.5\textwidth]{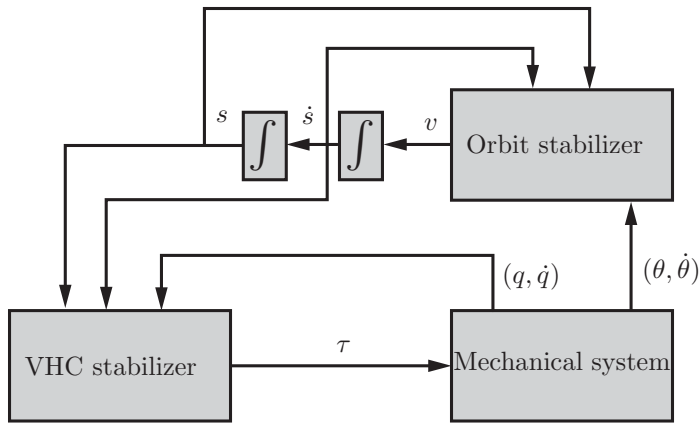}}
\caption{Block diagram of the \vhc-based orbital stabilizer.}
\label{fig:stab:s-Method}
\end{figure}

The block diagram of the \vhc-based orbital stabilizer is depicted in
Figure~\ref{fig:stab:s-Method}.
\begin{rem}
In Theorem~\ref{thm:main_result} we only claim asymptotic stability of
$\bar \gamma$, even though Theorem~\ref{thm:transverse_linearization}
guarantees that $\bar \gamma$ is {\em exponentially} stable relative
to $\bar \Gamma$. The reason is that in the proof we use a reduction
theorem for asymptotic stability of sets~\cite{el-hawwary-2013}. A
different proof technique could be used to show that $\bar \gamma$ is
in fact exponentially stable for the closed-loop system.  \hfill
$\triangle$
\end{rem}

\section{Discussion}\label{sec:stab:shir}

In this section we briefly compare the control methodology of this
paper with the ones
in~\cite{morris2009hybrid,canudas2002orbital,canudas2004concept,ShiPerWit05}.

\textbf{Comparison with~\cite{morris2009hybrid}.} The notion of
dynamic hybrid extension introduced by Morris and Grizzle
in~\cite{morris2009hybrid} bears a conceptual resemblance to dynamic
\vhcs and their extended reduced dynamics presented in this
article. In~\cite{morris2009hybrid}, the \vhcs that induce stable
walking gaits of biped robots are parameterized using variables whose
evolution are event-triggered. In particular, the \vhc parameters get
updated after each impact of the swing leg with the ground. The update
law is designed such that the invariance of a suitably modified
manifold, which the authors call the extended zero dynamics manifold,
is preserved while simultaneously enforcing a periodic stable walking
gait on the biped. Our approach follows the same philosophy of
preserving the invariance of a suitably modified manifold in order to
maintain the desired configurations of the mechanical system. However,
in our framework, the dynamics of the \vhc parameter are continuous
rather than event-triggered.

\textbf{Comparison with~\cite{canudas2002orbital}.} The approach by
Canudas-de-wit et al. in~\cite{canudas2002orbital} also relies on
dynamically changing the geometry of \vhcs. A target orbit on the
constraint manifold, which is generated by a harmonic oscillator, is
considered and the dynamics of the \vhc parameter is designed such
that the target orbit is stabilized on the constraint manifold. This
approach, however, cannot be used to stabilize an assigned closed
orbit induced by the original \vhc on the constraint
manifold. Moreover, the methodology in~\cite{canudas2002orbital} has
only been employed to control the periodic motions of a pendubot.  It
is unclear to what extent it can be generalized to other mechanical
systems.

\textbf{Comparison
  with~\cite{canudas2004concept,ShiPerWit05,ShiFreGus10}.}
In~\cite{canudas2004concept,ShiPerWit05}, the authors employ \vhcs to
find feasible closed orbits of underactuated mechanical systems. Once
the orbit is found, it is stabilized through transverse linearization
of the $2n$-dimensional dynamics~\eqref{eq:stab:ELsys} along the
closed orbit. Similarly to this paper,
in~\cite{canudas2004concept,ShiPerWit05} the stabilization of the
transverse linearization is carried out by solving a periodic Riccati
equation. But while the linearized system
in~\cite{canudas2004concept,ShiPerWit05} has dimension $2 n -1$, the
linearized system~\eqref{eq:transverse_linearization} always has
dimension $3$. And while the feedback in Theorem~\ref{thm:main_result}
is time-independent, the one proposed
in~\cite{canudas2004concept,ShiPerWit05} is
time-varying. Additionally, while the approach proposed in this paper
gives explicit parametrizations of the orbits to be stabilized, the
approaches in~\cite{canudas2004concept,ShiPerWit05} require the
knowledge of the actual periodic trajectory which is not available in
analytic form.  The most important difference between the approach in
this paper and the ones in~\cite{canudas2004concept,ShiPerWit05} lies
in the fact that, in~\cite{canudas2004concept,ShiPerWit05}, the
time-varying controller does not preserve the invariance of the
constraint manifold. The work in~\cite{ShiFreGus10} generalizes the
theory of~\cite{ShiPerWit05} to systems with degree of underactuation
greater than one. The philosophy in~\cite{ShiFreGus10} is analogous to
that of~\cite{ShiPerWit05} and shares the same differences just
outlined with our work. The authors use virtual constraints to help
identify a desired closed orbit of the control system, then linearize
the control system around said orbit to design a stabilizer. In this
paper, we only deal with systems with degree of underactuation one.

We end this section with a remark about the computational cost of the
controller proposed in Theorem~\ref{thm:main_result}. The controller
has two components: an input-output feedback linearizing controller,
$\tau^\star$, enforcing the dynamic \vhc, and a scalar feedback,
$v^\star$, for the double-integrator stabilizing the desired closed
orbit. The computational cost of these controllers for real-time
implementation is essentially equivalent to that of virtual constraint
controllers used by Grizzle and collaborators for biped robots and by
many other researchers in the area. The design of the orbit stabilizer
$v^\star$ involves the solution of a periodic Riccati equation for the
three-dimensional linear periodic
system~\eqref{eq:transverse_linearization}. The dimension of this
problem is always $3$, independent of the number of DOFs of the
original mechanical system. As described above, this is a major
advantage of the simultaneous stabilization method proposed in this
paper.  We surmise that the proposed approach can be particularly
effective to reduce the design complexity for robots with a large
number of DOFs.
\section{Example}\label{sec:MotEx}

\begin{figure}[ht]
	\begin{center}
\psfrag{x1}[.25cm]{$x [m]$} \psfrag{x3}[-0.25cm]{$y [m]$}
\includegraphics[scale=0.4]{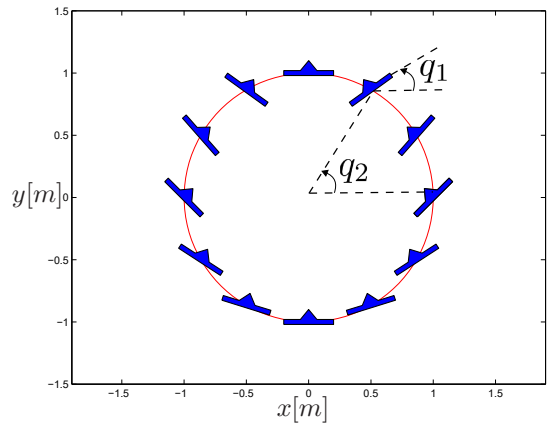}
	\end{center}
	\caption{Configurations of a PVTOL vehicle on the unit circle
          under the \vhc proposed in~\cite{consolini2010path}.}
	\label{fig:stab:PVTOL_circle}
\end{figure}

In this section we use the theory developed in this paper to enhance a
result found in~\cite{consolini2010path}. We consider the model of a
V/STOL aircraft in planar vertical take-off and landing mode (PVTOL),
introduced by Hauser et al. in~\cite{HauSasMey92}. The vehicle in
question is depicted in Figure~\ref{fig:stab:PVTOL_circle}, where it
is assumed that a preliminary feedback has been designed making the
centre of mass of the aircraft lie on a unit circle on the vertical
plane, $\mathcal{C}=\{x\in \Re^2 :\, \big|x\big|=1\}$, also depicted
in the figure. In~\cite{consolini2010path} it was shown that the model
of the aircraft on the circle is given by
\begin{equation}
\begin{aligned}
\ddot{q}_1 &= \frac{\mu}{\epsilon}
\big(g\sin(q_1)-\cos(q_1-q_2)\dot{q}_2^2+\sin(q_1-q_2)u\big),\\ \ddot{q}_2
&= u,
\end{aligned}
\end{equation}
where $q_1$ denotes the roll angle, $q_2$ the angular position of the
aircraft on the circle, and $u$ the so-called tangential control input
resulting from the design in~\cite{consolini2010path}.  Also, $\mu$
and $\epsilon$ are positive constants. In this example, we set $\mu /
\epsilon=1$.

\begin{figure}[ht]
\psfrag{a}[bc]{$\gamma^+$} \psfrag{b}[cc]{$\gamma^-$}
\psfrag{q}{$q_2$} \psfrag{Q}{$\dot q_2$}
\centerline{\includegraphics[scale=0.4]{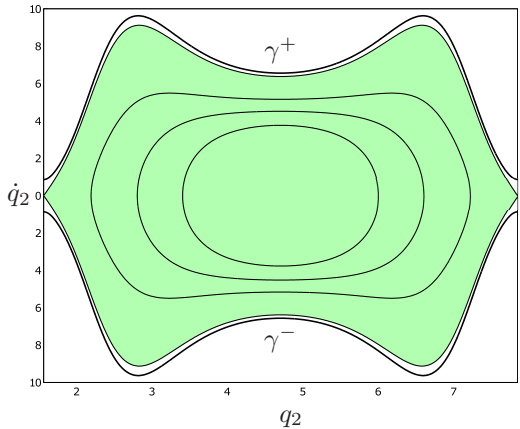}}
\caption{The phase portrait of the reduced dynamics of the PVTOL
  vehicle under the \vhc depicted in
  Figure~\ref{fig:stab:PVTOL_circle}. The closed orbits in the shaded
  area correspond to oscillations. The rest of the orbits correspond
  to rotations. We would like to stabilize the counterclockwise
  rotation $\gamma^{+}$ corresponding to the energy level set
  $E_0=41.5$.}
\label{fig:stab:PVTOL_phase}
\end{figure}

In~\cite{consolini2010path}, a feedback $u(q,\dot q)$ was designed to
enforce a regular \vhc of the form $h(q) = q_1 - f(q_2)=0$,
represented in Figure~\ref{fig:stab:PVTOL_circle}. It was shown that
the ensuing reduced dynamics, a few orbits of which are depicted in
Figure~\ref{fig:stab:PVTOL_phase}, are Lagrangian. Each closed orbit
in Figure~\ref{fig:stab:PVTOL_phase} represents a motion of the PVTOL
on the circle, with roll angle $q_1$ constrained to be a function of
the position, $q_2$, on the circle. Orbits in the shaded area
represent a rocking motion of the PVTOL along the circle (these are
oscillations), while orbits in the unshaded area represent full
traversal of the circle (these are rotations).  The theory
in~\cite{consolini2010path} was unable to stabilize individual closed
orbits of the reduced dynamics. The theory of this paper fills the gap
left open in~\cite{consolini2010path}.

We wish to stabilize the closed orbit $\gamma^{+}$ depicted in
Figure~\ref{fig:stab:PVTOL_phase} which corresponds to the energy
level set $E_0=41.5$. The parametrization of $\gamma^+$ on the
$(q_2,\dot q_2)$ plane is $\s \mapsto (\s,\sqrt{2/M(\s) (E_0 -
  V(\s))})$. Here, $T_1 = T_2 = 2\pi$.  We render the \vhc dynamic by
setting $h^s(q) = q_1-L_1s - f(q_2 - L_2s)=0$, with $L = \col(L_1,L_2)
= \col(1, 1)$. We enforce this dynamic \vhc by means of the
feedback $\tau^\star$ in~\eqref{eq:stab:tau_s_method}, with
$k_p=100$ and $k_d=10$. Thus trajectories converge to $\bar \Gamma$
at a rate of $\exp(-5t)$. This rate of convergence is chosen so as
to make the enforcement of the dynamic \vhc faster than the orbit
stabilization mechanism.

Since $\gamma^+$ is a rotation, we parameterize it with the
map~\eqref{eq:param_rotation}.  We check numerically that the
$2\pi$-periodic pair $(A(t),B(t))$
in~\eqref{eq:transverse_linearization},~\eqref{eq:A_and_b_components}
is controllable, and after some tuning we pick $R = 400$ and $Q =
\diag\{1/2,\, 10^4,\, 1\}$ to set up the Riccati
equation~\eqref{eq:PRE}.  We numerically solve this equation using the
one-shot generator method~\cite{hench1994numerical} (see
also~\cite{johansson2009tools,gusev2010numerical} for a detailed
treatment of existing numerical algorithms to solve the periodic
Riccati equation) and find the gain matrix $K(\cdot)$. The resulting
characteristic multipliers of the transverse
linearization~\eqref{eq:transverse_linearization:general} with
time-varying feedback $u=K(t) z$ are $\{0.0447, -3.6816 \times 10^{-5}
\pm 2.7122\times 10^{-5} i\}$. This means that trajectories on the
constraint manifold near $\bar \gamma$ converge to $\bar \gamma$ at
a rate of $\exp[\log(0.0447)t /(2 \pi) ] = \exp(-0.49 t)$. Thus the
enforcement of the dynamic \vhc occurs faster than the orbit
stabilization mechanism.

\begin{figure}[ht]
	\begin{center}
		\psfrag{t}{$t \ \text{[sec]}$}
		\psfrag{h}[-.3cm]{$h^{s(t)}(q(t))$}
		\includegraphics[scale=0.4]{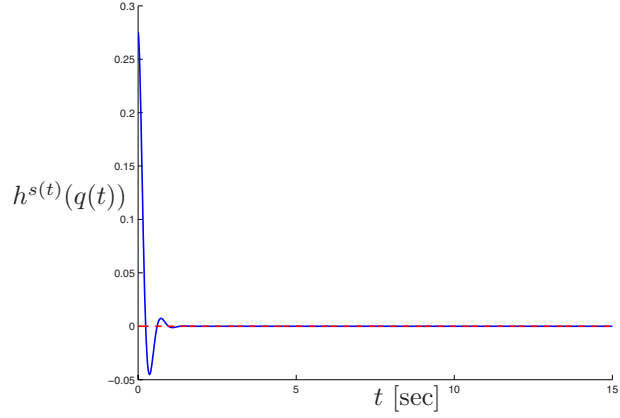}
	\end{center}
	\caption{The dynamic \vhc $h^s(q)=0$ is asymptotically stabilized on the vehicle.}
	\label{fig:PVTOL_vhcerr}
\end{figure}

\begin{figure}[ht]
	\begin{center}
		\psfrag{t}{$t \ \text{[sec]}$}
		\psfrag{s}{$s(t)$}
		\includegraphics[scale=0.4]{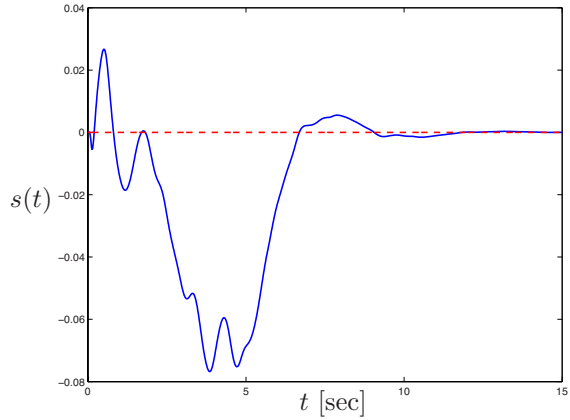}
	\end{center}
	\caption{Output of the double integrator.}
	\label{fig:fig_s}
\end{figure}

\begin{figure}[ht]
	\begin{center}
		\psfrag{t}[c]{$t \ \text{[sec]}$}
		\psfrag{E}[c]{$E(q(t),\dot q(t))$}
		\includegraphics[scale=0.4]{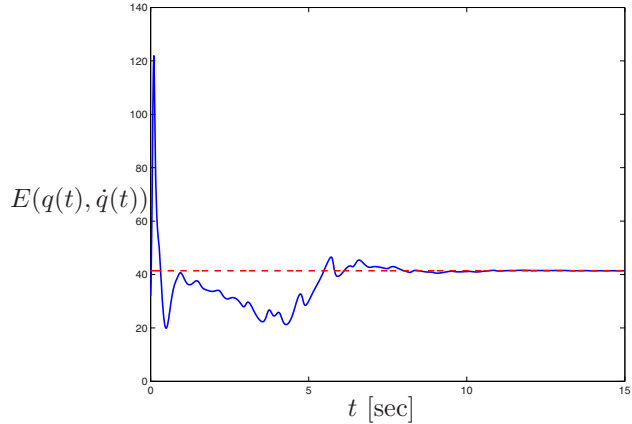}
	\end{center}
	\caption{Energy of the vehicle on the constraint manifold.}
	\label{fig:figE}
\end{figure}

\begin{figure}[ht]
	\begin{center}
                \psfrag{A}[r]{$2 \pi \sin(q_2)$}
                \psfrag{B}{$2 \pi \cos(q_2)$}
                \psfrag{C}{$\dot q_2$}
                \psfrag{G}{$\gamma^+$}
		\includegraphics[scale=0.4]{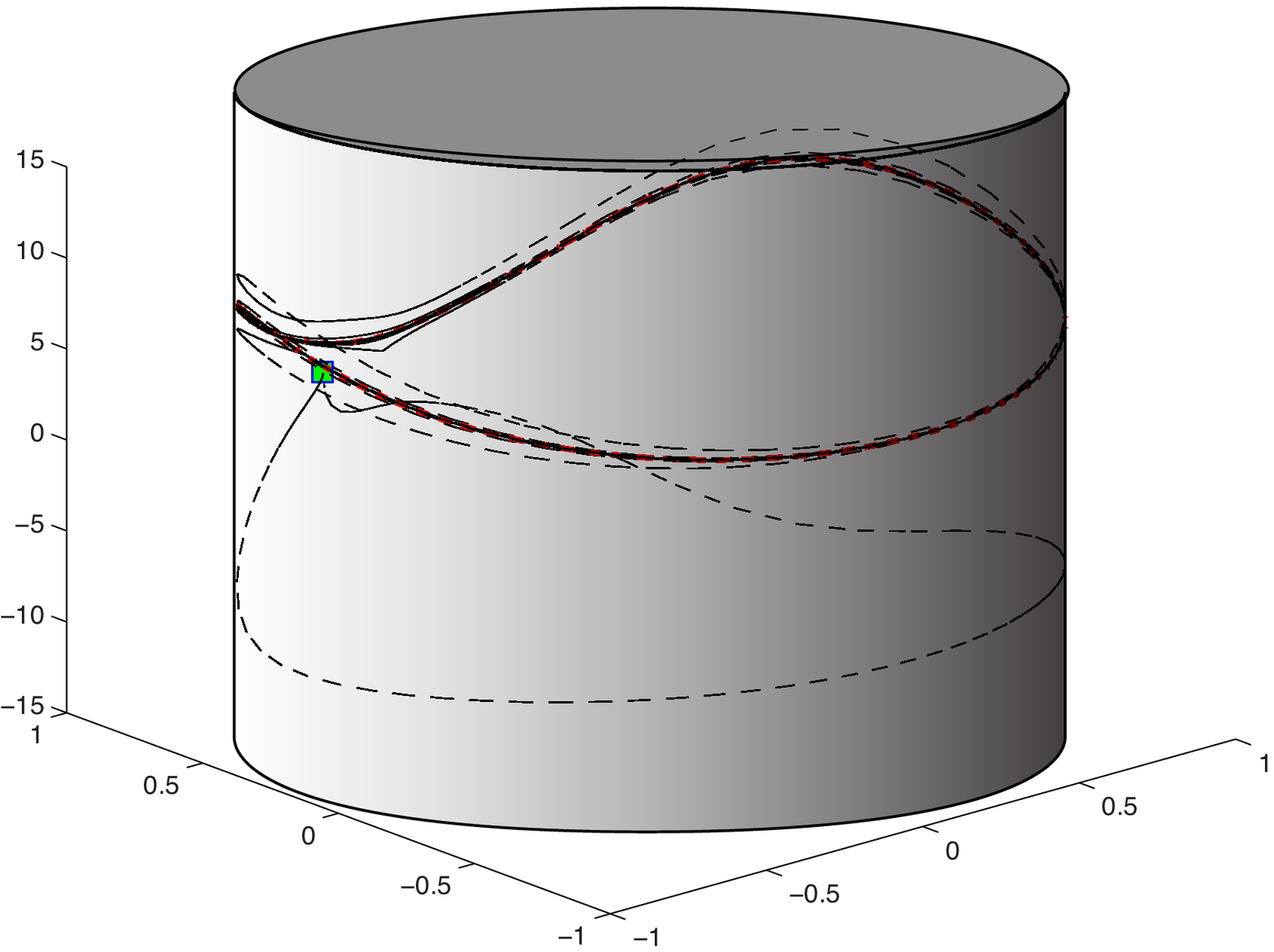}
	\end{center}
	\caption{The time trajectory of $(q_2,\dot{q}_2)$ on the 
		cylinder $\ucirc\times \Re$.}
	\label{fig:stab:PVTOL_cylinder}
\end{figure}

\begin{figure}[tbp]
	\begin{center}
	    \psfrag{t}{$t \ \text{[sec]}$}
            \psfrag{q}[c]{$q_1(t)$}	
            \includegraphics[scale=0.4]{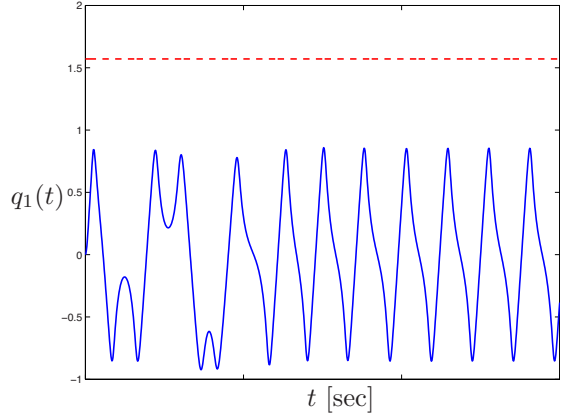}
	\end{center}
	\caption{The time trajectory of $q_1$.}
	\label{fig:stab:PVTOL_q1}
\end{figure}

The simulation results for the controller in
Theorem~\ref{thm:main_result} are presented next. We pick the initial
condition $q(0)=(0,\pi/2+0.2)$, $\dot q(0)=(0,0)$, $(s(0),\dot
s(0))=(0,0)$. We verified that other initial conditions in a
neighborhood of $\bar \Gamma$ give similar results as the ones that
follow.  Figures~\ref{fig:PVTOL_vhcerr} and~\ref{fig:fig_s} depict the
graph of the function $h^{s(t)}(q(t))$ and the output of the double
integrator, $s(t)$, respectively. They reveal that the \vhc is
properly enforced and that $s(t) \to 0$. Figures~\ref{fig:figE}
and~\ref{fig:stab:PVTOL_cylinder} depict the energy of the vehicle on
the constraint manifold and the time trajectory of
$(\theta(t),\dot{\theta}(t))$ on the cylinder $\ucirc\times \Re$. The
energy level $E_0$ is stabilized and the trajectory on the cylinder
converges to $\gamma^+$. Finally, Figure~\ref{fig:stab:PVTOL_q1}
depicts the graph of the roll angle $q_1(t)$, demonstrating that, due
to the enforcement of the dynamic version of the \vhc depicted in
Figure~\ref{fig:stab:PVTOL_circle}, the vehicle does not roll over
for the given initial condition. When comparing
Figure~\ref{fig:PVTOL_vhcerr} with Figures~\ref{fig:fig_s}
and~\ref{fig:figE}, it is evident that the enforcement of the
dynamic \vhc occurs faster than the convergence to the closed
orbit. As a final remark, in the proposed framework the roots
of the polynomial
$s^2 + k_d s + k_p$ determine the rate of convergence of trajectories to the
constraint manifold $\bar \Gamma$, while the characteristic
multipliers concern the constrained dynamics on $\bar \Gamma$, and
they characterize the rate at which trajectories on $\bar \Gamma$
converge to the closed orbit $\bar \gamma$.

\section{Conclusions}

We have proposed a technique to enforce a \vhc on a mechanical control
system and simultaneously stabilize a closed orbit on the constraint
manifold. The theory of this paper is applicable to mechanical control
systems with degree of underactuation one. For higher degrees of
underactuation, the reduced dynamics are described by a differential
equation of order higher than two and, generally, the problem of
characterizing closed orbits becomes harder. The result of
Section~\ref{sec:stab:s_control} concerning the exponential
stabilization of closed orbits for control-affine systems is still
applicable in this case.

\appendix

\section{Proofs of Technical Results in Section~\ref{sec:stab:dynConstr}}\label{sec:appendixA}

{\sc \bf Proof of Proposition~\ref{prop:stab:augmented}.}
Considering the output $e=h(q-Ls)$ and taking two derivatives along
system~\eqref{eq:stab:sys_augmented}, we get
\[
\ddot{e} = (\star) -dh\big|_{q-Ls}Lv + A^{s}(q)\tau,
\]
where $A^{s}(q)=dh_{q-Ls}D^{-1}(q)B(q)$. Denote $\mu^s:=\min_{q \in
  (h^s)^{-1}(0)} \det A^s(q)$. Then, $s\mapsto \mu^s$ is a continuous
function. We claim that $\mu^0 \neq 0$. Indeed, the assumption that
$h(q)=0$ is regular implies by Proposition~\ref{prop:regular_vhc} that
$\det A^0(q) \neq 0$ for all $q \in (h^0)^{-1}(0) = h^{-1}(0)$. Since
$h^{-1}(0)$ is a compact set and $q \mapsto \det A^0(q)$ is
continuous, $\min(\det A^0(q)) \neq 0$, proving that $\mu^0\neq 0$, as
claimed. By continuity, there exists an open interval $\cI \subset
\Re$ containing $s=0$ such that $\mu^s \neq 0$ on $\cI$ implying that
$A^{s}(q)$ is nonsingular for all $q \in (h^s)^{-1}(0)$ and all $s \in
\cI$.  \hfill\(\qed\)

{\sc \bf Proof of Proposition~\ref{prop:stab:stabilizability}.}
By~\eqref{eq:stab:Gammabar}, we have $\bar{\Gamma} = \{
(q,\dot{q},s,\dot{s})\in T\bar{\cQ} : (q-Ls,\dot{q}-L\dot{s})\in
\Gamma\}$, from which it follows that
$\|(q,\dot{q},s,\dot{s})\|_{\bar{\Gamma}} =
\|(q-Ls,\dot{q}-L\dot{s})\|_\Gamma$. This fact and the
inequalities in~\eqref{eq:classK} imply that
\begin{equation}\label{eq:stab:lemprop}
\alpha(\|(q,\dot{q},s,\dot{s})\|_{\bar{\Gamma}})\leq
H(q-Ls,\dot{q}-L\dot{s}) \leq
\beta(\|(q,\dot{q},s,\dot{s})\|_{\bar{\Gamma}}).
\end{equation}
Letting $e=h(q-Ls)$, the feedback~\eqref{eq:stab:tau_s_method} gives
$\ddot{e}+k_d\dot{e}+k_p e = 0$, so that the equilibrium
$(e,\dot{e})=(0,0)$ is asymptotically stable. Since $(e,\dot{e})=
H(q-Ls,\dot{q}-L\dot{s})$, property~\eqref{eq:stab:lemprop} implies
that $\bar{\Gamma}$ is asymptotically stable.\hfill\(\qed\)

\section{Proof of Theorem~\ref{thm:transverse_linearization}}\label{sec:appendixB}

Let $H:\cX \to \Re^{n-1}$ and $\pi: \cU \to [\Re]_T$ be as in the
theorem statement.  We claim that there exists a neighborhood $\cV$ of
$\gamma$ in $\cX$ such that the map $F: \cV\to [\Re]_T \times
\Re^{n-1}$, $x \mapsto (\s,\z)=(\pi(x),H(x))$ is a diffeomorphism onto
its image. By the generalized inverse function
theorem~\cite{GuiPol:74}, we need to show that $dF_x$ is an
isomorphism for each $x \in \gamma$, and that $F|_\gamma$ is a
diffeomorphism $\gamma \to [\Re]_T \times \{0\}$. The first property
was proved in~\cite[Proposition 1.2]{HauChu94}. For the second
property, we observe that $F|_\gamma = \pi|_\gamma \times \{0\}$ is a
diffeomorphism $\gamma \to [\Re]_T \times \{0\}$, since $\pi|_\gamma
=\varphi^{-1}$ is a diffeomorphism $\gamma \to [\Re]_T$. The smooth
inverse of $F|_\gamma$ is
\begin{equation}\label{eq:T_inverse}
(F|_\gamma )^{-1} = F^{-1}(\s,0) = \varphi(\s).
\end{equation}
Thus $F:\cV \to [\Re]_T \times \Re^{n-1}$ is a diffeomorphism onto its
image, as claimed.  Since $\s \mapsto \varphi(\s)$ is a regular
parameterization of the orbit $\gamma$, and since $\gamma$ is an
invariant set for the open-loop system, $f(\varphi(\s))$ is
proportional to $\varphi'(\s)$. More precisely, defining the continuous
function $[\Re]_T \to \Re$,
\begin{equation}\label{eq:rho}
\rho(\s)=\frac{\langle f(\varphi(\s)),\varphi'(\s) \rangle} {\|\varphi'(\s)\|^2},
\end{equation}
we have that
\begin{equation}\label{eq:proportionality_term}
(\forall \s \in [\Re]_T) \ \varphi'(\s) = \frac{ 1}{\rho (\s)}f(\varphi(\s)),
\end{equation}
and $\rho$ is bounded away from zero. We now represent the control
system~\eqref{eq:control_affine} in $(\s,\z)$ coordinates. The
development is a slight variation of the one presented in the proof
of~\cite[Proposition 1.4]{HauChu94}, the variation being due to the
fact that, in~\cite{HauChu94}, it is assumed that $\rho=1$. For the
$\s$-dynamics, we have
\[
\dot \s = \big[ L_f \pi(x) + L_g \pi(x) u \big]_{x=F^{-1}(\s,\z)}.
\]
We claim that the restriction of the drift term to $\gamma$ is
$\rho(\s)$. Indeed, using~\eqref{eq:T_inverse}
and~\eqref{eq:proportionality_term}, we have
\[
\big[ L_f \pi(x) \big]_{x=F^{-1}(\s,0)} = L_f \pi(\varphi(\s)) =
d\pi_{\varphi(\s)} f(\varphi(\s)) = \rho(\s) d\pi_{\varphi(\s)}
\varphi'(\s) = \rho(\s).
\]
The last equality is due to the fact that $\pi(\varphi(\s))=\s$,
so that $d\pi_{\varphi(\s)} \varphi'(\s)=1$. Thus we may write
\[
\dot \s = \rho(\s) + f_1(\s,\z) + g_1(\s,\z) u,
\]
where $f_1(\s,0)=0$. The derivation of the $\z$ dynamics is
essentially the same as in~\cite[Proposition 1.4]{HauChu94} so we
present their form without proof. The control
system~\eqref{eq:control_affine} in $(\s,\z)$ coordinates has the form
\begin{equation}\label{eq:transformed_dyn}
\begin{aligned}
& \dot \s = \rho(\s) + f_1(\s,\z) + g_1(\s,\z) u\\
& \dot \z = \bar A(\s)\z+f_2(\s,\z) + g_2(\s,\z)u,
\end{aligned}
\end{equation}
where $f_1$ and $f_2$ satisfy $f_1(\s,0)=0$, $f_2(\s,0)=0$,
$\partial_\z f_2(\s,0)=0$.

Letting $\T = \int_0^T \big| 1/\rho(u) \big| du$, we have that $\T >0$
because $\rho$ is bounded away from zero. Consider the partial
coordinate transformation $\tau: [\Re]_T \to [\Re]_{\T}$ defined as
\[
\tau(\s) = \left[ \int_0^\s 1/ \rho(u) du \right]_{\T}.
\]
Since $\rho$ is bounded away from zero, the derivative $\tau'(\s)$ is
also bounded away from zero, implying that $\tau$ is a
diffeomorphism. We denote by $\s(\tau)$ the inverse of
$\tau(\s)$. System~\eqref{eq:transformed_dyn} in $(\tau,\z)$
coordinates reads as
\begin{equation}\label{eq:transformed_dyn2}
\begin{aligned}
& \dot \tau = 1 + \tilde f_1(\tau,\z) + \tilde g_1(\tau,\z) u\\
& \dot \z = \bar A(\s(\tau))\z+ \tilde f_2(\tau,\z) +
  g_2(\s(\tau),\z)u,
\end{aligned}
\end{equation}
where $\tilde f_1(\tau,\z) = f_1(\s(\tau),\z) / \rho(\s(\tau))$,
$\tilde g_1(\tau,\z) = g_1(\s(\tau),\z)/\rho(\s(\tau))$, and $\tilde
f_2(\tau,\z) = f_2(\s(\tau),\z)$.

System~\eqref{eq:transformed_dyn2} has the same form of that
in~\cite[Proposition 1.4]{HauChu94} (which, however, has no control
inputs). By~\cite[Proposition 1.5]{HauChu94}, we deduce that the orbit
$\gamma$ is exponentially stabilizable if and only if the
$\T$-periodic system
\[
\frac{dz}{d\tau} = \bar A(\s(\tau)) z + \tilde g_2(\s(\tau),0) u,
\]
is stabilizable. Since $\s(\tau)$ is a diffeomorphism, we may perform
the time-scaling
\begin{equation}\label{eq:transverse_linearization:general:proof}
\frac{dz}{d \s} = \frac{1}{\rho(\s)} \left[ \bar A(\s) z + g_2(\s,0) u\right].
\end{equation}
Thus $\gamma$ is exponentially stabilizable if and only if the
$T$-periodic system~\eqref{eq:transverse_linearization:general:proof}
is asymptotically stable. By comparing the system and input matrices
of~\eqref{eq:transverse_linearization:general:proof} with those of
system~\eqref{eq:transverse_linearization:general}, we see that to
prove part (a) of Theorem~\ref{thm:transverse_linearization} it
suffices to show that
\begin{align} \label{eq:A(t)}
\bar A(\s) & =[(dL_f H)_{\varphi(\s)}] dH^\dagger_{\varphi(\s)}\\
g_2(\s,0) & =  L_g H(\varphi(\s)) \label{eq:B(t)}.
\end{align}
Since $\z = H(x)$, the coefficient of $u$ in $\dot \z$ is 
\[
g_2(\s,\z) =  L_g H \circ F^{-1}(\s,\z).
\]
Using~\eqref{eq:T_inverse} we get $g_2(\s,0)= L_g H \circ
F^{-1}(\s,0) = L_g H (\varphi(\s))$. This proves
identity~\eqref{eq:B(t)}.

Concerning identity~\eqref{eq:A(t)}, and referring to
system~\eqref{eq:transformed_dyn}, $\bar A(\s)$ is the Jacobian of
$\dot \z$ with respect to $\z$ evaluated at $(\z,u)=(0,0)$. Since
\[
\dot \z = L_f H \circ F^{-1}(\s,\z) + L_g H \circ
F^{-1}(\s,\z) u,
\]
we have
\[
\bar A(\s) = \partial_\z \big[ L_f H \circ F^{-1}(\s,\z)\big]\big|_{\z=0}.
\]
By the chain rule and the identity~\eqref{eq:T_inverse}, we get
\[
\bar A(\s) = [(d L_f H)_{\varphi(\s)}] \partial_\z F^{-1}(\s,\z) \big|_{\z=0}.
\]
To show that identity~\eqref{eq:A(t)} holds, we need to show that
$\partial_\z F^{-1}(\s,\z) \big|_{\z=0} =
dH_{\varphi(\s)}^\dagger$. To this end, we use the fact that
\[
dF_{\varphi(\s)} dF^{-1}_{(\s,0)}  = I_{n},
\]
or
\[
\begin{bmatrix} d\pi_{\varphi(\s)} \\ 
dH_{\varphi(\s)}
\end{bmatrix} [\partial_\s F^{-1} 
  \ \ \partial_\z F^{-1}(\s,\z) ]\Big|_{\z=0} =I_n.
\]
In light of the above, $\partial_\z F^{-1}(\s,\z)
\big|_{\z=0}$ is uniquely defined by the identities
\[
\begin{aligned}
& d\pi_{\varphi(\s)} \partial_\z F^{-1}(\s,\z) \big|_{\z=0} =0
  \\
& dH_{\varphi(\s)} \partial_\z F^{-1}(\s,\z)
  \big|_{\z=0} = I_{n-1},
\end{aligned}
\]
so we need to show that 
\begin{align}
&  d\pi_{\varphi(\s)} dH_{\varphi(\s)}^\dagger =0 \label{eq:first}\\
& dH_{\varphi(\s)} dH_{\varphi(\s)}^\dagger = I_{n-1}, \label{eq:second}
\end{align}
Identity~\eqref{eq:second} holds by virtue of the fact that
$dH^\dagger$ is the right-inverse of $dH$. Using the definition of
pseudoinverse and taking the transpose of~\eqref{eq:first}, we may
rewrite~\eqref{eq:first} as
\[
dH_{\varphi(\s)} d\pi_{\varphi(\s)}\trans =0.
\]
Since $\varphi(\pi(x))=x$ for all $x \in \gamma$, we have $d\varphi_{\pi(x)} d\pi_x = I_n$, or
\[
(\forall x \in \gamma) \ d\pi_x\trans = \frac{d \varphi_{\pi(x)}}{\|d\varphi_{\pi(x)}\|_2^2},
\]
so that 
\[
dH_{\varphi(\s)} d\pi_{\varphi(\s)}\trans =
\frac{dH_{\varphi(\s)}d
  \varphi_{\s}}{\|d\varphi_{\s}\|_2^2}.
\]
Since $H(\varphi(\s)) \equiv 0$, $dH_{\varphi(\s)} d
\varphi_\s =0$ for all $\s \in [\Re]_T$. Thus,
$dH_{\varphi(\s)} d\pi_{\varphi(\s)}\trans =0$ for all $\s
\in [\Re]_T$. 

We have thus shown that identities~\eqref{eq:first}
and~\eqref{eq:second} hold, implying that identity~\eqref{eq:A(t)}
holds. This concludes the proof of part (a) of the theorem.

For part (b), let $A(t)$, $B(t)$ be as
in~\eqref{eq:transverse_linearization:general}, and suppose that the
origin of $\dot \z =(A(t) + B(t)K(t)) \z$ is asymptotically
stable. With the controller $u^\star(x) = K(\pi(x))H(x)$, the dynamics
of the closed-loop system in $(\tau,\z)$ coordinates read as
\begin{equation}\label{eq:transformed_coordinates_cls}
\begin{aligned}
& \dot \tau = 1 + \tilde f_1(\tau,\z) + \tilde g_1(\tau,\z) K(\s(\tau)) \z \\
& \dot \z = \bar A(\s(\tau)) +\tilde f_2(\tau,\z) + g_2(\s(\tau),\z) K(\s(\tau)) \z.
\end{aligned}
\end{equation}
For the $\z$ dynamics we have
\begin{align}
\dot \z &= [\bar A(\s(\tau)) + g_2(\s(\tau),0)K(\s(\tau))]\z + \tilde f_2(\tau,\z) +
     [g_2(\s(\tau),\z) - g_2(\s(\tau),0)]K(\s(\tau))z \nonumber \\
&=  [\bar A(\s(\tau)) + g_2(\s(\tau),0)K(\s(\tau))]\z + \tilde F_2(\s,\z), 
\label{eq:z_dynamics}
\end{align}
with $\tilde F_2(\s,0)=0$, $\partial_\z \tilde F_2(\s,0)=0$. By using
$\s$ as time variable, the linear part of the $z$-dynamics reads as
\[
\frac{dz}{d\s} = \frac{1}{\rho(\s)} [\bar A(\s) + g_2(\s,0)K(\s)]\z.
\]
Using the identities~\eqref{eq:A(t)} and~\eqref{eq:B(t)} we rewrite the above as
\[
\frac{dz}{d\s} = \big( A(\s)+B(\s) K(\s) \big) z.
\]
By assumption, the origin of this system is asymptotically stable,
implying that the origin of the system $\dot z =[\bar A(\s(\tau)) +
  g_2(\s(\tau),0)K(\s(\tau))]\z$ has the same property.  Referring
to~\eqref{eq:z_dynamics} and using~\cite[Proposition~1.5]{HauChu94},
we conclude that the closed orbit $\gamma$ is exponentially stable for
the closed-loop system~\eqref{eq:transformed_coordinates_cls} and
hence also for system~\eqref{eq:control_affine} with feedback
$u^\star(x)=K(\pi(x)) H(x)$. \hfill\(\qed\)

\begin{ack}                               
A. Mohammadi and M. Maggiore were supported by the
Natural Sciences and Engineering Research Council (NSERC) of
Canada. A. Mohammadi was partially supported by the University of Toronto 
Doctoral Completion Award (DCA).  
\end{ack}
\bibliographystyle{plain}       
{\small \bibliography{PhDBiblio}}
\end{document}